\documentclass[letterpaper, 10 pt, conference]{ieeeconf}

\IEEEoverridecommandlockouts                            
\usepackage{graphics} 
\usepackage{epsfig} 
\usepackage{times} 
\usepackage{amsmath} 
\usepackage{amssymb}  
\newtheorem{defi}{\textbf{Definition}}
\newtheorem{remark}{Remark}

\newtheorem{assume}{\textbf{Assumption}}
\newtheorem{lemma}{Lemma}
\usepackage[ruled,vlined]{algorithm2e}
\usepackage[utf8]{inputenc}
\usepackage{cite}
\usepackage{amsmath,amssymb,amsfonts}
\usepackage{algorithmic}
\usepackage{graphicx}
\usepackage{textcomp}
\usepackage{xcolor}
\usepackage{comment}
\usepackage{tikz}
\usepackage{accents}
\usepackage{multirow}
\newcommand*{\dt}[1]{%
   \accentset{\mbox{\large\bfseries .}}{#1}}

\def\BibTeX{{\rm B\kern-.05em{\sc i\kern-.025em b}\kern-.08em
    T\kern-.1667em\lower.7ex\hbox{E}\kern-.125emX}}

\title{\LARGE \bf 
Unified Control Framework: A Novel Perspective on Constrained Optimization, Optimization-based Control, and Parameter Estimation}

\author{Revati Gunjal, Syed Shadab Nayyer, Sushama Wagh, and Navdeep Singh 
\thanks{Revati Gunjal and Syed Shadab Nayyer are joint first authors. Revati Gunjal, Sushama Wagh, and Navdeep Singh are with the Control and Decision Research Center (CDRC), EED, Veermata Jijabai Technological Institute, Mumbai 400019, India. Syed Shadab Nayyer is working as an Electrification--Specialist in the Transportation Business Unit (TBU) at TATA Elxsi, India.
        {\tt\small rgunjal\_p21@ee.vjti.ac.in, \tt\small shadgec@gmail.com} ** This is a preliminary version, the revised manuscript will be released soon.}%
 }

\begin{document}

\maketitle

\begin{abstract}
A common theme in all the above areas is designing a dynamical system to accomplish desired objectives, possibly in some predefined optimal way. Since control theory advances the idea of suitably modifying the behavior of a dynamical system, this paper explores the role of control theory in designing efficient algorithms (or dynamical systems) related to problems surrounding the optimization framework, including constrained optimization, optimization-based control, and parameter estimation. This amalgamation of control theory with the above mentioned areas has been made possible by the recently introduced paradigm of Passivity and Immersion (P\&I) based control. The generality and working of P\&I, as compared to the existing approaches in control theory, are best introduced through the example presented below.

\end{abstract}
\begin{keywords}
 Constrained Optimization, Control Barrier Function, Gradient Descent, Optimization-based control, Parameter Estimation, Primal-Dual formulation, Passivity and Immersion (P\&I) approach
\end {keywords}
\section{Introduction}
Optimization plays a significant role in various domains, including machine learning, data science, game theory, and a wide range of scientific and engineering applications. Optimization problems are extensively dealt with using the Gradient Descent (GD) method. According to the requirements, various variants of the GD method are proposed in the literature. These methods are usually executed in the form of discrete algorithms. When these algorithms are perceived in continuous-time systems, they are equivalent to the numerical integration of the differential equations corresponding to a dynamical system \cite{gunjal2024nesterov}. The dynamical perspective of the algorithms aids in the modification of these algorithms to ensure the improvement in the associated convergence characteristics. Since the themes (constrained optimization, optimization-based control, and parameter estimation) discussed in the paper can be visualized as dynamical systems, performance enhancement through the control input targeting a particular objective seems a natural and intuitive choice.

The general approach of the paper is to geometrically visualize the particular optimization problem as a Manifold Stabilization Problem (MSP). Manifold stabilization involves designing a control scheme to guarantee the convergence (attractivity) of the system dynamics to a specific invariant manifold \cite{gunjal2023new}. The manifold is selected targeting a specific objective that can be seeking the optimal for optimization, or equilibrium stabilization within a feasible region of operation as required in the optimization-based control, or convergence of parameter estimates to the actual system parameters as desired in the parameter estimation problem. Following the selection of an invariant manifold, its attractivity is achieved through a control scheme formulated using the recently introduced Passivity and Immersion (P\&I) approach \cite{nayyer2022towards} \cite{nayyer2022passivitylcss}. P\&I is a systematic and efficient technique to solve a specific problem by visualizing it as a manifold stabilization problem. The embedding of the geometry through an invariant manifold makes the technique generalized enough to accommodate a wide variety of problems in systems and control including stabilization, control \cite{nayyer2022passivitylcss}, synchronization \cite{nayyer2024synchronization}, optimization \cite{gunjal2024controlled} \cite{gunjal2024nesterov}, parameter estimation \cite{nayyer2022passivityPE}, to name a few. 

The geometry of the loss landscape plays an important role in ensuring the convergence to the minima in the optimization problems. While dealing with the optima-seeking problems on the curved space, the knowledge of the geodesics (distance between the trajectories lying on the curved surface) plays a vital role. The calculation of geodesics demands information about the Riemannian metric that describes the system geometry. If a state-independent Riemannian metric exists, then the geodesics are just straight lines \cite{boyd1994linear}. However, finding a state-dependent Riemannian metric is a difficult problem. The literature mentions various approaches for the evaluation of the Riemannian metric.  The Control Contraction Metric (CCM) approach proposed by Ian Manchester et al. \cite{manchester2017control} evaluates the Riemannian metric (geodesics) by finding the minimal length path that joins the current state to the desired state. This real-time optimization problem is solved with the Sum-of-Squares (SOS) approach and requires solvers such as Mosek and the parser YALMIP \cite{lopez2019contraction}. Another popular approach proposed in \cite{wang2015immersion} \cite{wang2016immersion} explores the use of Horizontal Finsler-Lyapunov functions \cite{forni2013differential} that decay along the trajectories of the prolonged system to ensure the contraction between system trajectories. The P\&I approach depends on the pseudo-Riemannian metric (PR) \cite{gunjal2024controlled}, which is a generalization of a Riemannian metric with a relaxed requirement of positive definiteness. The splitting of the tangent bundle with the PR metric makes the P\&I approach a unified framework to accommodate a variety of problems.

The most popular choice of algorithm for solving constrained optimization problems is the continuous-time primal-dual-gradient dynamics (PDGD) technique \cite{kose1956solutions}. Due to the straightforward formulation and scalability, PDGD has been widely adopted in fields like power systems \cite{zhao2014design} \cite{li2015connecting}, wireless communication \cite{chen2011convergence} \cite{chiang2007layering}, distributed resource allocation \cite{ding2018primal} \cite{ferragut2013network} and image processing \cite{chambolle2011first}. Theoretic analysis of the performance of PDGD, particularly the convergence property, recently gained considerable attention. In the case of an unconstrained optimization problem where the objective function to optimize is strongly convex and smooth, the exponential stability of the gradient dynamics is ensured. Global exponential stability is an appealing aspect in practice. Strong stability guarantees are desired for the control systems in critical infrastructures like the power grid. While solving the optimization algorithms, discretization is an essential step. The global exponential stability ensures that the simple explicit Euler discretization has a geometric convergence rate when the discretization step size is sufficiently small \cite{Stuart1994NumericalAO} \cite{stetter1973analysis}.

The PDGD of convex optimization with the strongly convex and smooth objective function with non-strongly convex constraints is considered in this paper. The unconstrained problem with strongly convex and smooth objectives always attains an exponential convergence to the optimal points. In the case of constrained optimization, the Lagrangian function constructed with the strongly convex objectives and concave constraints is not strongly convex and hence does not converge exponentially to the equilibrium. The variants of PDGD with the various modifications in the Lagrangian function such as the augmentation as proposed in  \cite{qu2018exponential}, by adding a regularization term to the Lagrangian function as proposed in \cite{bianchin2021time}, or by multiplying the PDGD dynamics with a Riemannian metric inspired by the natural gradient formulation as proposed in \cite{bansode2019exponential}. Although these variants of PDGD succeeded in achieving exponential convergence, the optimal solution obtained with each of them is approximate due to the modifications in the Lagrangian function. 

To achieve the exponential convergence to the optimal solution without any approximation, the control theoretic perspective (Controlled PGGD) is proposed in the paper. The paper emphasizes the analysis of continuous-time gradient-based optimization through a control-theoretic perspective via a passivity and immersion (P\&I) approach. In the case of a time-invariant objective, the fundamental analysis of gradient descent based on convexity and the natural gradient is replaced by a more general analysis through a control perspective. Due to the notion of an invariant manifold and its attractivity, the proposed Controlled PDGD becomes a generalized technique to accommodate a variety of constrained optimization problems with convex and non-convex objective functions with linear, non-linear, convex, and nonconvex constraints. The projection operator on the feasible set is used to ensure the convergence of optimal solution within the feasible region. The projection operator is insufficient when the constraints are non-convex. In such scenarios, the proposed Controlled PDGD method is extended by modifying the invariant manifold such that the requirement of the projection operator is avoided. This approach inspires the control perspective in optimization-based control. 

Optimization-based control (Safety-critical control) problems often involve a strong interconnection between the likely conflicting control objectives and safety constraints. Such safety-critical control problems are often dealt with using the Quadratic Programming (QP) framework \cite{ames2012control} mediating between the conflicting objectives i.e. stability and safety \cite{ames2014rapidly}. In \cite{ames2016control} the Control Lyapunov Function (CLF) (defined to express the control objective) and Control Barrier Function (CBF) (defined to express the safety objective) are unified through the QP. The QP mediates the trade-off of achieving stabilization while remaining within the feasible (safe) set using relaxation. In this formulation, finding the CLF and CBF is a tedious procedure, and solving the QPs increases the computation complexity. Hence, a more natural and systematic control perspective is proposed in the paper. Since the problem is in the form of constrained optimization using a PDGD method is an intuitive choice rather than using a QP formulation. In the proposed approach the stability and safety objectives are represented in the form of PDGD dynamics. The stability ensuring the safety is achieved through the CBF-inspired invariant manifold and its attractivity. Through the notion of tangent bundle splitting with PR metric, a non-quadratic manifold control Lyapunov function is formulated. The control law evaluated through this Lyapunov function prevents the optimizer from leaving the safe region of operation.   

The third major theme of the paper is the Parameter Estimation (PE) problem, which plays a vital role in the development of adaptive control strategies and system identification.  The parameter estimation problem often deals with finding the unknown parameters of a plant (or controller) with the knowledge of the measurable signals. Classically, the PE problem is solved using the Gradient Descent technique. The parametric error equation (PEE) is formulated using the estimated parameters with the parameter estimator and the actual parameters of the system. The optimization problem is solved to find exact estimates of the parameters such that the PEE reduces to zero. This optimization problem leads to a linear time-varying (LTV) system. To ensure the global exponential convergence of PEE of such systems the regressor vector must satisfy the persistency of excitation (PE) condition \cite{ortega2020modified}. Satisfying the PE property is very rare, hence to alleviate this challenge a relaxed and weaker notion of Interval excitation (IE) is used in the Concurrent Learning-based estimators \cite{chowdhary2010concurrent}. 

For achieving faster convergence to the actual parameters a filtering approach is proposed alternatively in the Dynamic Regressor Extension and
Mixing (DREM) \cite{ortega2020new} procedure. The implementation of a DREM algorithm for $`\mathrm{n}'$ unknown parameters requires '$\mathrm{n-1}$' number of filters that puts forward challenges like poor tuning of filter coefficients and learning rate \cite{korotina2020parameter}. To avoid the complications of implementing multiple filters and the choice of higher learning rates, a control perspective is proposed in the paper for designing the parameter estimator. The gradient estimator is cast as a Manifold stabilization problem to design a novel framework of Controlled Gradient Estimator (CGE). The proposed approach is computationally advantageous as it uses only one filter to construct a memory regressor extension (MRE) \cite{kreisselmeier1977adaptive}. Due to the attractivity of the manifold, faster convergence of the estimated parameters to the actual parameters is achieved, which further relieves the necessity of high learning rates.

The remaining paper is structured as follows: The control perspective for constrained optimization (Controlled-PDGD) along with a wide variety of numerical examples is demonstrated in Section \ref{Control in PD}. 
The control perspective for the optimization-based control with the formulation of CBF-inspired PDGD dynamics and the numerical example is illustrated in Section \ref{control in CBF}. Section \ref{Control in PE} represents the formulation of the Controlled Gradient Estimator (CGE) along with the representative examples. The paper is concluded in Section \ref{conclusion}. The step-wise procedure of the Passivity and Immersion (P\&I) approach is mentioned in the Appendix.

\section{A Control Perspective in Constrained Optimization} \label{Control in PD}
\subsection{Controlled PDGD dynamics}
Consider a constrained optimization problem
\begin{equation}\label{optimization problem}
    \begin{split}
  &\mathrm{min}\; \mathrm{f}(\mathrm{x})\\
  &\mathrm{s.t.}\; \mathrm{x}\in \mathrm{X}\\
  &\mathrm{g_i}(\mathrm{x})\leq 0,\forall ^\mathrm{m}_{\mathrm{i}=1}
    \end{split}
\end{equation}
where 
\begin{equation}
    \mathrm{X}=\left \{ \mathrm{x}\in \mathbb{R}^\mathrm{n}|\mathrm{g_i}(\mathrm{x})\leq 0,\forall ^\mathrm{m}_{\mathrm{i}=1} \right \}
\end{equation}
is the feasible set of $\mathrm{x}$ (a decision variable).
The functions $\mathrm{f}:\mathbb{R}^\mathrm{n}\rightarrow \mathbb{R}$ (objective function) and $\mathrm{g}:\mathbb{R}^\mathrm{n}\rightarrow \mathbb{R}^\mathrm{m}$ (a function vector describing the inequality constraints) are assumed continuously differentiable $\left ( \mathit{C}^2 \right )$ with respect of $\mathrm{x}$, such that the following assumptions hold.
\begin{assume} \label{A1}
The gradient of function $\mathrm{f}$, $\nabla\mathrm{f}:\mathbb{R}^\mathrm{n}\rightarrow \mathbb{R}^\mathrm{n}$ is strongly monotone on the domain $\mathrm{X}$, i.e., 
\begin{equation*}
    (\mathrm{x}_1-\mathrm{x}_2)^\mathrm{T}(\nabla\mathrm{f}(\mathrm{x}_1)-\nabla\mathrm{f}(\mathrm{x}_2))\geq \mu \left \| \mathrm{x}_1-\mathrm{x}_2 \right \|^2
\end{equation*}
This implies that the objective function $\mathrm{f}$ is strongly convex in $\mathrm{x}$ with the modulus of convexity, $\frac{\mu}{2}$, while $\mu> 0$.
\end{assume}


The above assumption ensures that $\mathrm{x}$ is strictly feasible and strong duality holds for the optimization problem (\ref{optimization problem}).

\begin{defi} \label{strong convexity}
\textit{(Strong convexity). A twice differentiable function $\mathrm{f:\mathbb{R}^n\rightarrow \mathbb{R}}$ is $\mu$-strongly convex with $\mu>0$ if its Hessian matrix $\mathrm{\nabla^2f(x)}$ satisfies the matrix inequality,
\begin{equation}
\mathrm{\nabla^2f(x)\geq \mu I , \hspace{1cm} \mathrm{x\in \mathbb{R}^n} }    
\end{equation}}
\end{defi}

The Lagrangian $\mathrm{L}(\mathrm{x},\lambda ):\mathbb{R}^\mathrm{n}\times \mathbb{R}^\mathrm{m}\rightarrow \mathbb{R}$ associated with the optimization problem (\ref{optimization problem}) is given as 
\begin{equation} \label{Lagrangian}
    \mathrm{L}(\mathrm{x},\lambda )=\mathrm{f(x)+\lambda ^Tg(x)}
\end{equation}
where $\lambda \in \Lambda \subseteq \mathbb{R}^\mathrm{m}_+=\left \{ \lambda \in \mathbb{R}^\mathrm{m},\lambda_i\geq 0,\forall ^\mathrm{m}_{\mathrm{i}=1}  \right \}$ represents the vector of Lagrangian multipliers associated with the constraint functions $\mathrm{g(x)}$. The above defined Lagrangian $\mathrm{L(x,\lambda)}$ is $\mathbb{C}^2$- differentiable convex-concave in $\mathrm{x}$ and $\lambda$ respectively, i.e., $\mathrm{L(.,\lambda)}$ is convex for all $\lambda \in \Lambda$ and $\mathrm{L(x,.)}$ is concave for all $\mathrm{x\in X}$. The point $\mathrm{(x^*, \lambda^*)}$ is a saddle-point, if it follows the condition,
\begin{equation}
    \mathrm{L(x^*,\lambda)\leq L(x^*, \lambda^*)\leq L(x,\lambda^*)}\;\;\; \forall \mathrm{x\in X},\; \forall \lambda\in \Lambda
\end{equation}

Let $\mathrm{z=(x,\lambda)\in \Omega =X\times \Lambda }$, where $\Omega$ is a non-empty and closed subset of $\mathbb{R}^\mathrm{n}$ and $\mathbb{R}^\mathrm{m}_{\geq 0}$. $\mathrm{z^*=(x^*,\lambda^*)}$ is the saddle point solution of (\ref{Lagrangian}).

The Primal-Dual Gradient Dynamics (PDGD) associated with the Lagrangian $\mathrm{L}(\mathrm{x},\lambda )$ mentioned above is given as
\begin{equation}
    \begin{split}
 \dot{\mathrm{x}}=-\nabla_{\mathrm{x}}\mathrm{L}(\mathrm{x},\lambda )=-\nabla \mathrm{f(x)}-\nabla (\mathrm{g(x))^T}\lambda\\\dot{\lambda}=\nabla_{\lambda}\mathrm{L}(\mathrm{x},\lambda )=\mathrm{g(x)}
    \end{split}
\end{equation}
\textcolor{blue}{The control perspective is explored by adding $\mathrm{u}$ in the $\dot{\lambda}$ i.e. the dynamics of the dual variable as a significant contribution of this paper.} The idea of stabilization and control is utilized in this optimization problem by externally adding the virtual control term $\mathrm{u}$. 
The PDGD with the addition of $\mathrm{u}$ is given by
\begin{equation}\label{PDGD}
    \begin{split}
\dot{\mathrm{x}}=-\nabla \mathrm{f(x)}-\nabla \mathrm{^Tg(x)}\lambda\\\dot{\lambda}=\mathrm{g(x)+u}
    \end{split}
\end{equation}

This input $\mathrm{u}$ is evaluated through the procedure of passivity and immersion (P\&I) approach \cite{nayyer2022towards} (The step-wise procedure is overviewed in Appendix) such that the solution of problem (\ref{optimization problem}) will converge exponentially to equilibrium. With the choice of the suitable implicit manifold, the considered constrained optimization problem is geometrically visualized as a Manifold Stabilization Problem (MSP). The controlled PDGD is formulated through the four steps of the MSP:\\
$\mathbf{(\mathcal{S}1)}$ \textbf{Target Dynamics :} Target dynamics for any system is the dynamics replicating a certain desired behavior of the system. The target dynamics might be asymptotically stable or exponentially stable or stable. In the considered optimization problem, the exponential behavior of the system trajectories is desired, hence the target dynamics is defined as $\dot{\mathrm{x}}=-\nabla\mathrm{f(x)}+\mathrm{k\nabla^Tg(x)g(x)}$. \\
$\mathbf{(\mathcal{S}2)}$ \textbf{Implicit Manifold :} The optimal solution to the problem under consideration should satisfy the conditions imposed on dual variables and constraints,
\begin{equation}\label{dual constraints}
    \lambda_i\geq 0,\forall ^\mathrm{m}_{\mathrm{i}=1}
\end{equation}
\begin{equation}\label{constraints condition}
    \mathrm{g_i}(\mathrm{x})\leq 0,\forall ^\mathrm{m}_{\mathrm{i}=1}
\end{equation}
Hence, the manifold should be selected such that the above conditions are satisfied. 
\begin{equation}\label{implicit}
    \lambda=- \mathrm{g}(\mathrm{x})
\end{equation}
Referring (\ref{implicit}), the implicit manifold is defined as 
\begin{equation}\label{manifold}
\Psi(\mathrm{x},\lambda)=\lambda+\mathrm{k}\mathrm{g}(\mathrm{x})=0.
\end{equation}
The implicit manifold will ensure that the convergence will be exponential with the flexible parameter $\mathrm{k}$.
\begin{remark}
    With the flexible parameter $\mathrm{k}$, the target dynamics can be defined as $\dot{\mathrm{x}}=-\nabla\mathrm{f(x)}+\mathrm{k\nabla^Tg(x)g(x)}$. \textcolor{blue}{With the smaller values of the parameter, $\mathrm{k}$ the contribution of the gradient of constraint function ($\mathrm{\nabla^Tg(x)g(x)}$) becomes negligible and the better minima is obtained.} 
\end{remark}
$\mathbf{(\mathcal{S}3)}$ \textbf{Invariance of the Manifold :} Conceptually, manifold invariance implies that the target dynamics should always lie on the implicit manifold. 
\begin{lemma} \label{Lemma 1}
    The implicit manifold $\Psi(\mathrm{x},\lambda)$ is invariant.
\end{lemma}
 \textit{Proof:} According to the relationship 
$\lambda=-\mathrm{kg(x)}$, the manifold is obtained as $\Psi(\mathrm{x},\lambda)=\lambda+\mathrm{k}\mathrm{g}(\mathrm{x})=0$ . The normal to the manifold is given as
\begin{equation} \label{normal}
  \nabla\Psi(\mathrm{x},\lambda)=\begin{bmatrix}
\mathrm{k\nabla g(x)} & \mathrm{I}
\end{bmatrix}
\end{equation}
The product of (\ref{normal}) and the velocity vector field $\begin{bmatrix}
\dot{\mathrm{x}}\\ 
\dot{\lambda}
\end{bmatrix}$ yields
\begin{equation}
   \begin{bmatrix}
\mathrm{k\nabla g(x)} & \mathrm{I} 
\end{bmatrix}\begin{bmatrix}
\dot{\mathrm{x}}\\ 
\dot{\lambda}
\end{bmatrix}=\begin{bmatrix}
\mathrm{k\nabla g(x)} & \mathrm{I} 
\end{bmatrix}\begin{bmatrix}
\dot{\mathrm{x}}\\ 
\mathrm{-k\nabla g(x) \dot{x}}
\end{bmatrix}=0
\end{equation}
 This implies that the velocity vector field is always tangent to the manifold $\Psi$. Hence, the implicit manifold $\Psi$ is invariant.\\
$\mathbf{(\mathcal{S}4)}$ \textbf{Attractivity of the Manifold through the Passivity and Immersion (P\&I) Approach :} After the construction of the invariant manifold, the convergence of the trajectories of off-the-manifold dynamics (manifold attractivity) is assured through the control input. The control input is related to a storage function, which in turn is obtained from the splitting (connection) of the tangent space in the fiber bundle, defining the given dynamical system using the tangent kernel (TK) (which is a degenerate two-form resembling a pseudo-Riemannian (PR) metric).
(For detailed procedure of P\&I approach refer to \cite{nayyer2022passivitylcss},\cite{nayyer2022towards}).

The (non-quadratic) candidate Lyapunov function (or manifold control Lyapunov function) $\mathbb{S}(\mathrm{x}, \lambda)$ (i.e., storage function) is defined as $\mathbb{S}(\mathrm{x}, \lambda)=\frac{1}{2}( \lambda+\mathrm{k}\mathrm{g}(\mathrm{x}))^2$. The convergence of the trajectories of the off-the-manifold dynamics to the implicit manifold at an exponential rate $\alpha>0$ is accomplished by selecting the condition $\dt{\mathbb{S}}\leq -\alpha\mathbb{S}$, which further provides the desired control law as 
\begin{equation}
      \mathrm{u}=\mathrm{-g(x)-k\nabla g(x)\dot{x}-\frac{\alpha}{2}(\lambda+kg(x))}\\
\end{equation}
\begin{small}
 \begin{equation}
    \Rightarrow \mathrm{u}=\mathrm{-g(x)-k\nabla g(x)(-\nabla \mathrm{f(x)}-\nabla \mathrm{^Tg(x)}\lambda)-\frac{\alpha}{2}(\lambda+kg(x))}
\end{equation}   
\end{small}
Hence, with the above-evaluated control law the controlled PDGD dynamics in (\ref{PDGD}) is defined as 
\begin{small}
  \begin{equation}\label{PDGD_with_control}
    \begin{split}
\dot{\mathrm{x}}=-\nabla \mathrm{f(x)}-\nabla \mathrm{^Tg(x)}\lambda\\
\dot{\lambda}=\mathrm{-k\nabla g(x)(-\nabla \mathrm{f(x)}-\nabla \mathrm{^Tg(x)}\lambda)-\frac{\alpha}{2}(\lambda+kg(x))}
    \end{split}
\end{equation}  
\end{small}

\subsection{Stability Analysis}
The GD dynamics for the primal variable is given as
\begin{equation}
    \dot{\mathrm{x}}=-\nabla \mathrm{f(x)}-\mathrm{\nabla^T}\mathrm{g(x)}\lambda
\end{equation}
The target dynamics defined on the implicit manifold $\lambda+\mathrm{kg(x)}=0$ is given as
\begin{equation} \label{Target dynamics}
    \dot{\mathrm{x}}=-\nabla\mathrm{f(x)}+\mathrm{k\nabla^T g(x)g(x)}
\end{equation}
\begin{lemma} The gradient dynamics (\ref{Target dynamics}) is exponentially stable. 
\end{lemma}
\textit{Proof :} Consider the Krasovskii-type Lyapunov candidate function $\mathrm{V(x)=\frac{1}{2}\dot{x}^T\dot{x}}$ with the assumption that the constraint function $\mathrm{g(x)=Ax-b}$ is linear in $\mathrm{x}$. 

Differentiating the above Lyapunov function along the trajectories of (\ref{Target dynamics}) gives
\begin{small}
    \begin{equation}
   \mathrm{\dot{V}=\dot{x}^T\left (-\nabla^2f(x)\dot{x}+k\nabla^2 g(x)g(x)\dot{x}  +k\nabla^T g(x)\nabla g(x)\dot{x}\right )}
\end{equation}
\end{small}

\begin{equation}
    \mathrm{\dot{V}\leq \dot{x}^T\left ( -\nabla^2 f(x)+k\nabla ^T g(x) \nabla g(x) \right )\dot{x}}
\end{equation}
Since the constraint function is assumed as $\mathrm{g(x)=Ax-b}$, then $\mathrm{\nabla^Tg(x)\nabla g(x)=A^TA}$ and $\mathrm{\nabla^2g(x)}=0$.
\begin{equation}
     \mathrm{\dot{V}\leq \dot{x}^T\left ( -\nabla^2 f(x)+kA^TA\right )\dot{x}}
\end{equation}
As per definition \ref{strong convexity}, the  Hessian matrix $\mathrm{\nabla^2f(x)}$ of strongly convex function $\mathrm{f}$ satisfies the matrix inequality, $\mathrm{\nabla^2f(x)\geq \mu I}$ with $\mu>0$. 
 
Also, the symmetric matrix $\mathrm{A^TA}$ follows the inequality, $\mathrm{A^TA\leq q_2I}$, where $\mathrm{q_2}$ is the largest eigen value of the matrix $\mathrm{A^TA}$ \cite{bansode2019exponential}.
\begin{equation}
    \mathrm{\dot{V}\leq \dot{x}^T\left ( -\mu I+kq_2I \right )\dot{x}}
\end{equation}

\begin{equation} \label{target exponential}
    \mathrm{\dot{V}\leq -\gamma V(x)}
\end{equation}
where $\gamma =\mathrm{\nabla^2f(x)-k\nabla^T g(x)\nabla g(x)}$ and $\gamma>0$.

The dynamics (\ref{Target dynamics}) is exponentially stable if
\begin{equation}
    \mathrm{\mu-kq_2}>0 
\end{equation}
where $\mathrm{k}$ is chosen as $\mathrm{k<<q_2}$.
The global attractivity of the implicit manifold would ensure that the proposed dynamics is globally exponentially stable. 

\subsection{Projected Primal-Dual Dynamics}
To ensure that the optimal solution to the problem (\ref{PDGD_with_control}) remains in the feasible region $ \mathrm{X}=\left \{ \mathrm{x}\in \mathbb{R}^\mathrm{n}|\mathrm{g_i}(\mathrm{x})\leq 0,\forall ^\mathrm{m}_{\mathrm{i}=1} \right \}$, the projection method introduced in \cite{bansode2019exponential} is used. 
Hence, the projection is applied to the primal variable dynamics as shown:
\begin{equation}\label{projection operator}
    \mathrm{\dot{x}=\beta\left \{ P_X \left ( x-\alpha_x\frac{\partial L(x,\lambda)}{\partial x} \right )-x \right \}}
\end{equation}
where, $\alpha_x$, $\beta>0$ are parameters adjusted to control stability and ensure convergence. For the sake of simplicity, we set $\alpha_x=\beta=0$. $\mathrm{P_X}$ is a minimum norm projection operator of the form
\begin{equation}
    \mathrm{P_X=\underset{y\in X}{argmin}}\left \| \mathrm{y-x} \right \|
\end{equation}
\begin{remark}
    The projection in (\ref{projection operator}) is Lipschitz continuous; unlike the other types of discontinuous projections mentioned in the literature \cite{chen2020model} \cite{nagurney2012projected}\cite{zhu2018projected}, which project the dynamics onto the tangent cone of the feasible set, and thus they need the sophisticated analysis tools for discontinuous dynamical systems. 
\end{remark}
The projection method used above is referred to as global projection \cite{chen2021model} (that is, with the initial condition $\mathrm{x(0)\in X}$, the trajectories of $\mathrm{x(t)\in X}$ for all $\mathrm{t\geq 0}$). The intuitive notion behind such type of projection is that (\ref{projection operator}) attempts to take a step forward with step-size $\alpha_x$ along the gradient descent direction, then checks whether the arrival point $\mathrm{x-\alpha_x\frac{\partial L(x,\lambda)}{\partial x}}$ is feasible to $\mathrm{X}$. If feasible, (\ref{projection operator}) reduces to the ordinary gradient descent dynamics $\mathrm{\dot{x}=-\alpha_x\frac{\partial L(x,\lambda)}{\partial x}}$, otherwise a projection is performed to guarantee the feasibility of the solution $\mathrm{x}$. 

Once the feasibility of the solution $\mathrm{x}$ is ensured using the projection method, then the feasibility of the dual variable $\lambda$ is ensured through the relationship defined by the invariant manifold (\ref{manifold}). 

\subsection{Numerical Example: Controlled PDGD}
To illustrate the versatility of the proposed approach, both convex and nonconvex problems with linear, nonlinear, convex, and nonconvex constraints are considered.\\

\subsubsection{\textbf{Strongly Convex Objective Function with Inequality Constraints}}
To gain a deeper insight into the formulation of controlled PDGD, a quadratic cost function $\mathrm{f(x)=x^TWx}$ with the affine inequality constraints $\mathrm{Ax\leq b}$ has been considered. Here, $\mathrm{W}$ and $\mathrm{A}$ are $\mathrm{n\times n}$ and $mathrm{m \times n}$  Gaussian random matrices, and $\mathrm{b}$ is a Gaussian random vector. For understanding the formulation $\mathrm{n=3}$ and $\mathrm{m=2}$ are considered.

The PDGD dynamics for the considered example is
\begin{equation} 
    \begin{split}
 \dot{\mathrm{x}}=-\nabla \mathrm{f(x)}-\nabla=\mathrm{-Wx-A^T\lambda} \\\dot{\lambda}=\mathrm{g(x)=\mathrm{Ax-b}}
    \end{split}
\end{equation}

The corresponding Controlled PDGD dynamics is
\begin{equation} \label{CPDGD1}
    \begin{split}
 \dot{\mathrm{x}}=-\nabla \mathrm{f(x)}-\nabla=\mathrm{-Wx-A^T\lambda} \\\dot{\lambda}=\mathrm{g(x)=\mathrm{Ax-b+u}}
    \end{split}
\end{equation}
where $\mathrm{u}\in \mathbb{R}^\mathrm{m}$. The invariant manifold for this problem is 
\begin{equation}\label{Manifold1}
\Psi(\mathrm{x},\lambda)=\lambda+\mathrm{k}\mathrm{g}(\mathrm{x})=\lambda+\mathrm{k(Ax-b)}=0.
\end{equation}
The above-defined manifold is invariant as depicted by Lemma 1. Since the target dynamics defined on the manifold is exponential,
the corresponding invariant manifold is exponentially stable and the trajectories of the dynamics defined on the manifold converge exponentially to the equilibrium point. For the convergence of the trajectories of the dynamics of the ambient (or off-the-manifold) space, the control input $\mathrm{u}$ is designed. The attractivity of the internally stable invariant manifold is achieved as per the (S4) through the P\&I approach.

The storage function or candidate Lyapunov function is obtained from the splitting (connection) of the tangent space in the fiber
bundle representing the given dynamics.

\begin{equation}
   \mathrm{S(x,\lambda)=\frac{1}{2}(\lambda+kAx-kb)^2}
\end{equation}
The exponential convergence of the trajectories of the off-the-manifold (or ambient space) dynamics to the invariant manifold at
the rate $\alpha>0$ is ensured with the condition
\begin{equation}
    \mathrm{\dot{S}\leq -\alpha S}
\end{equation}
With the above condition, the final control law is evaluated as
\begin{small}
   \begin{equation} \label{Control1}
\mathrm{u=-Ax-b+kAWx+kAA^T\lambda-\frac{\alpha}{2}(\lambda+kAx-kb)}
\end{equation} 
\end{small}
The simulation results of the above example with different cases are illustrated further.\\

{\textbf{Case 1: Optimal Solution $(\mathrm{x}^*)$ is present within the feasible set}}

\begin{figure}[ht] 
 \centering
     \includegraphics[ width=\linewidth]{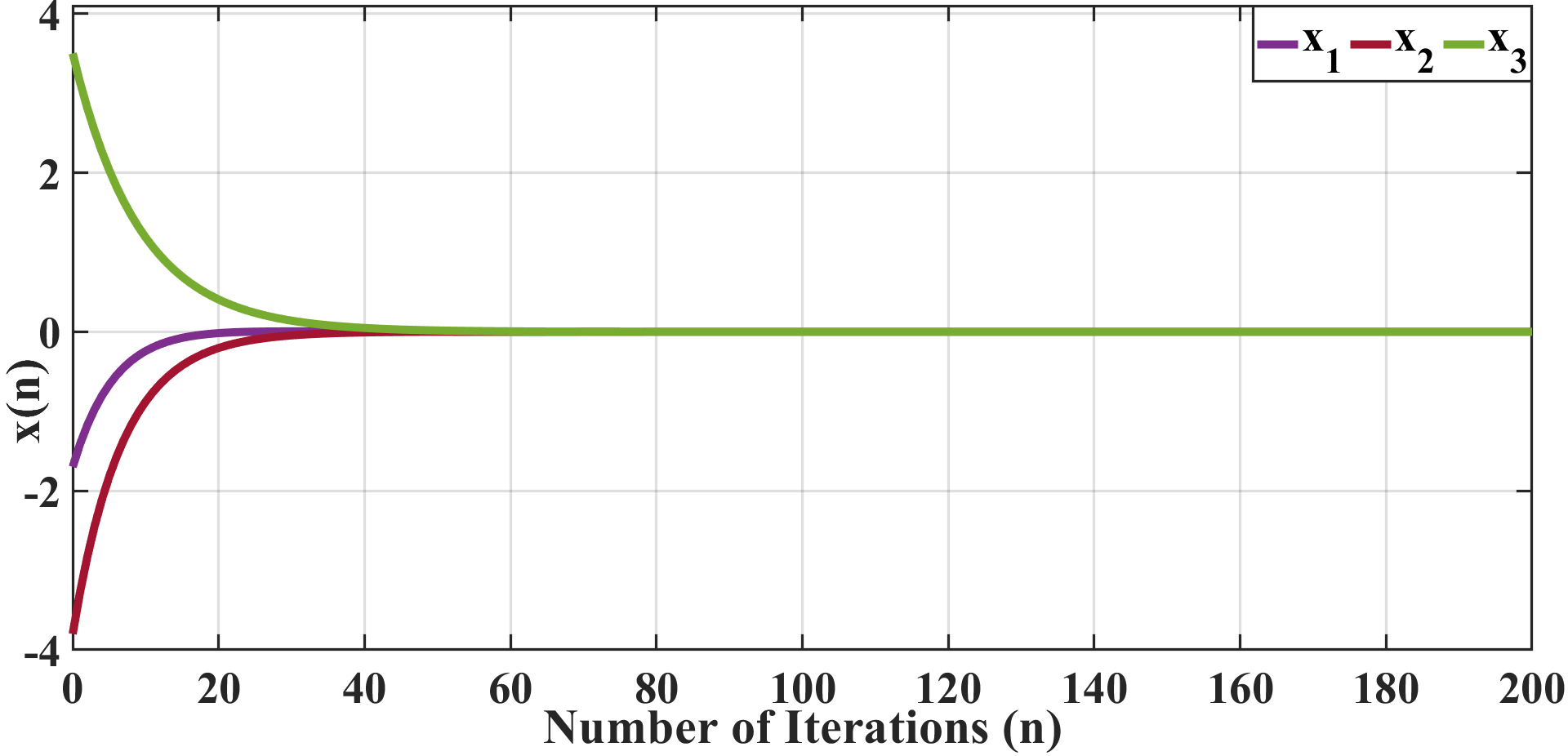}
     \caption{Evolution of the system dynamics in (\ref{CPDGD1}) with the control input (\ref{Control1}) for $\alpha=10$, $k=0.001$ and step size $\beta=0.01$. The y-axis illustrates the values of the variables, and the x-axis illustrates the number of iterations.}
     \label{xeqmin_x0in}
 \end{figure}

\begin{figure}[ht] 
 \centering
     \includegraphics[ width=\linewidth]{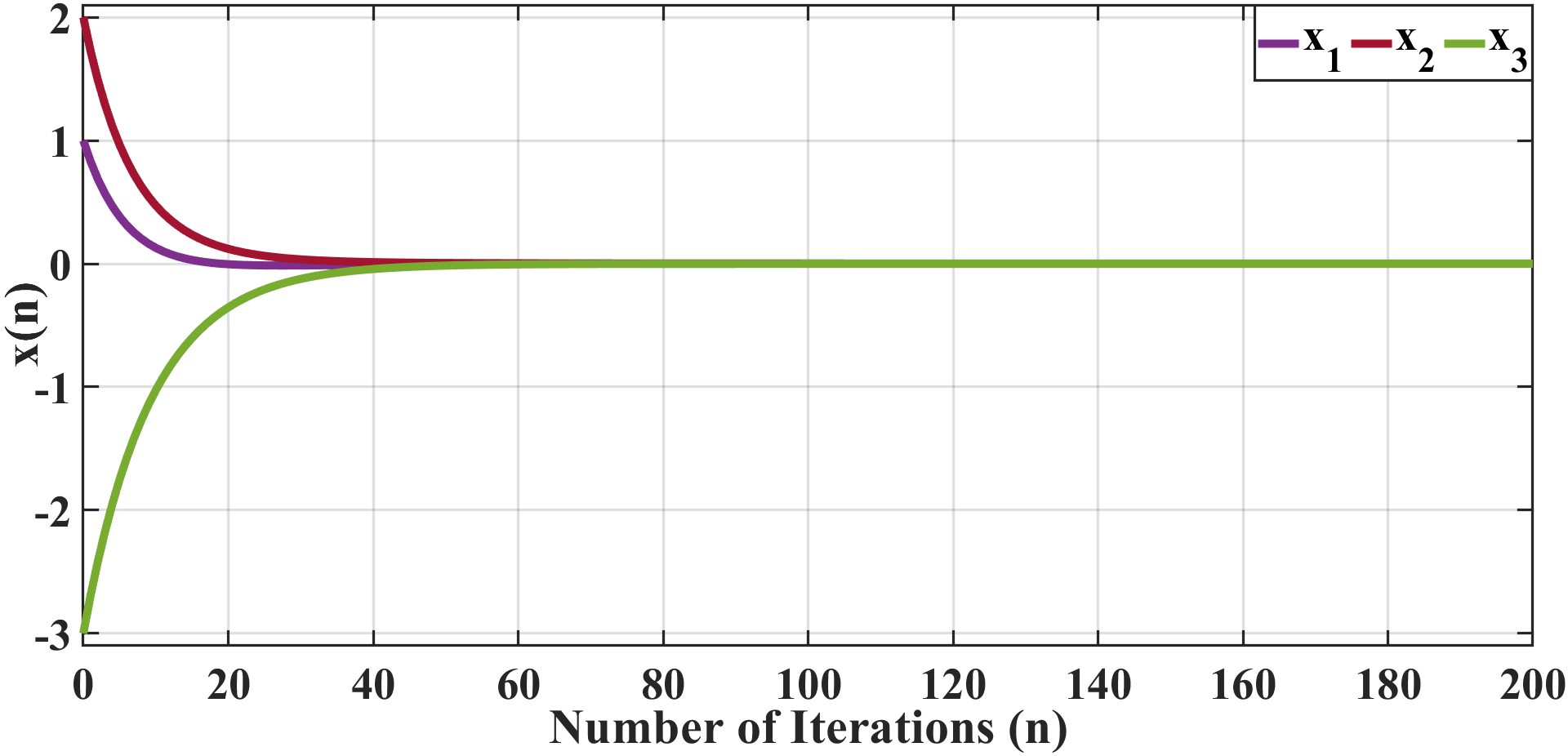}
     \caption{Evolution of the system dynamics in (\ref{CPDGD1}) with the control input (\ref{Control1}) for $\alpha=10$, $k=0.001$ and step size $\beta=0.01$. The y-axis illustrates the values of variables and the x-axis illustrates the number of iterations.}
     \label{xeqmin_x0out}
 \end{figure}

The time response of the system dynamics in (\ref{CPDGD1}) with the control law (\ref{Control1}) with the step size $\beta=0.01$ has been represented in the Fig.\ref{xeqmin_x0in} and Fig.\ref{xeqmin_x0out}. By varying the value of parameter $\alpha$ the convergence rate of the system trajectories can be adjusted. With the increasing value of $\alpha$, the rate of exponential convergence of the off-the-manifold trajectories increases. With the minimum value of the parameter $\mathrm{k}$ the exponential behavior of the primal variables is ensured while fulfilling the constraints. 

For the considered case, the constraints are selected such as the equilibrium point $\mathrm{x^*}$ lies within the feasible set. The projection operator is not necessary in such scenarios, since the optimizer will not leave the feasible set even if the initial point lies outside the feasible region as observed in Fig. \ref{xeqmin_x0out}.\\
 
{\textbf{Case 2: Optimal Solution $(\mathrm{x}^*)$ is present outside the feasible set: Necessity of projection}}

\begin{figure}[ht] 
 \centering
     \includegraphics[ width=\linewidth]{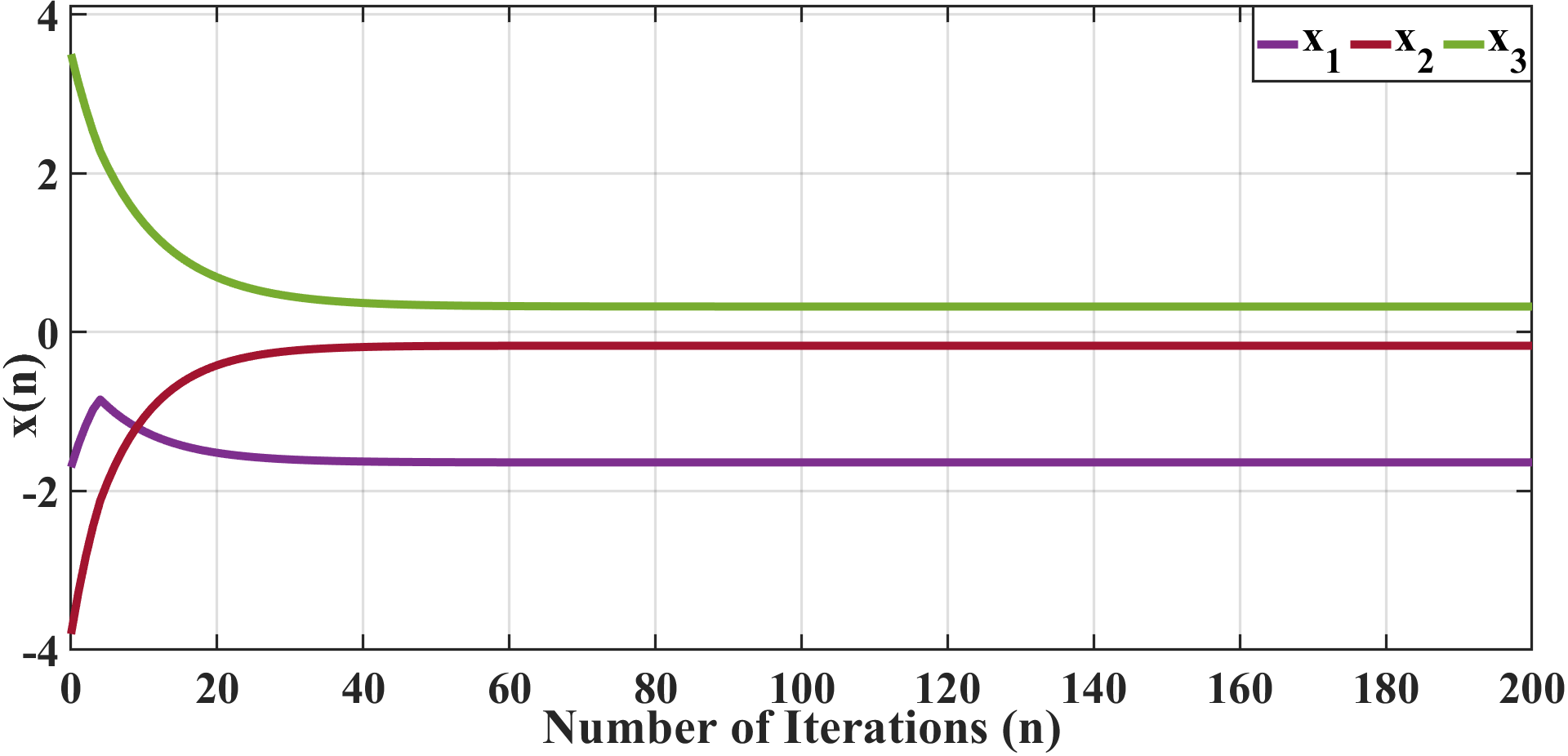}
     \caption{Evolution of the system dynamics in (\ref{CPDGD1}) with the control input (\ref{Control1}) for $\alpha=10$, $k=0.001$ and step size $\beta=0.01$. The y-axis illustrates the values of variables and the x-axis illustrates the number of iterations.}
     \label{xeqmout_x0in}
 \end{figure}

\begin{figure}[ht] 
 \centering
     \includegraphics[ width=\linewidth]{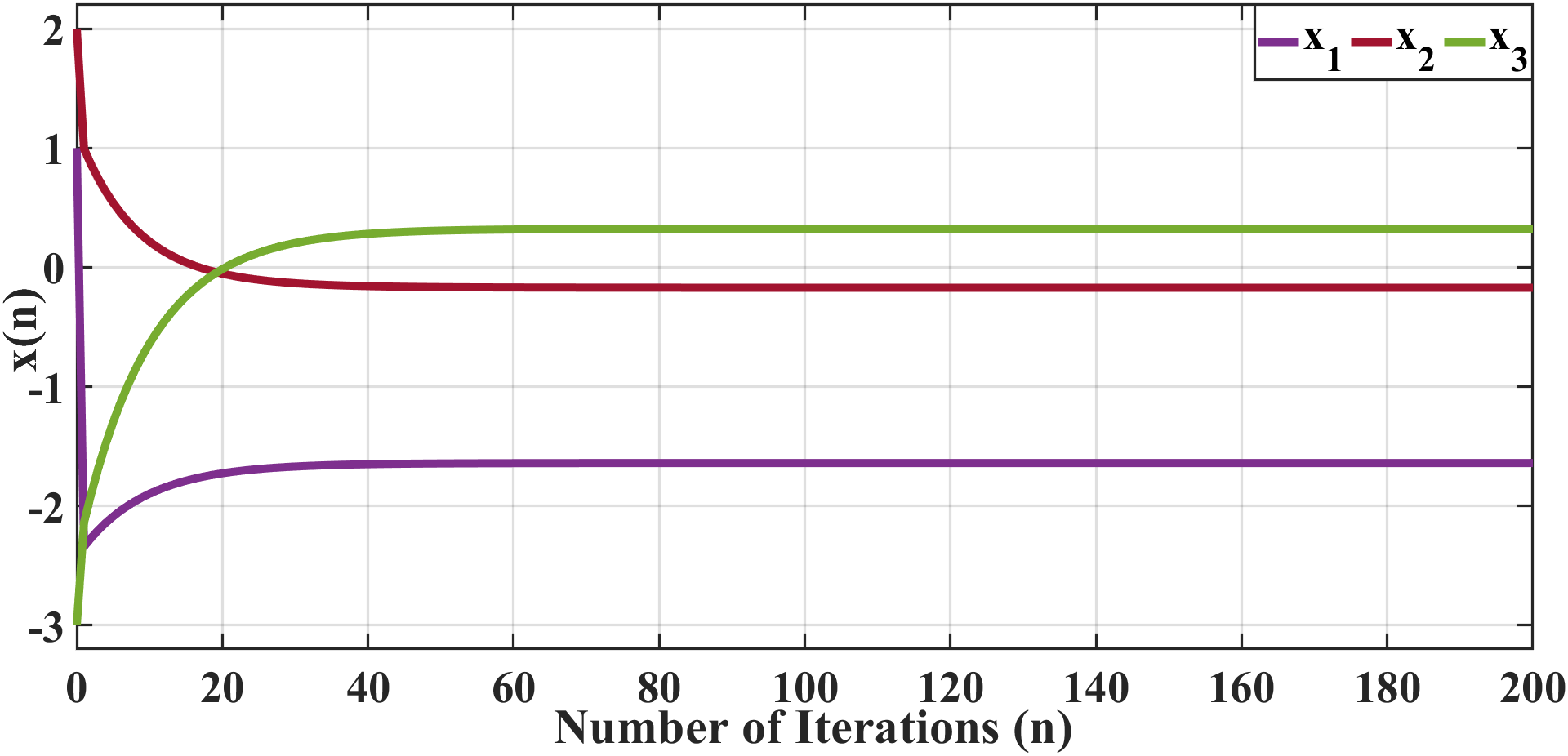}
     \caption{Evolution of the system dynamics in (\ref{CPDGD1}) with the control input (\ref{Control1}) for $\alpha=10$, $k=0.001$ and step size $\beta=0.01$. The y-axis illustrates the values of variables and the x-axis illustrates the number of iterations.}
     \label{xeqmout_x0out}
 \end{figure}

The constraints are selected such that the equilibrium point of the objective function lies outside the feasible region. In such scenarios, the optimizer tends to move outside the feasible region while searching for an optimal point. Hence, a projection operator is required to ensure that the constraints are satisfied. From Fig. \ref{xeqmout_x0in} and Fig. \ref{xeqmout_x0out}, it is evident that the optimizer is restrained from leaving the feasible region. Since we are using the projection on the feasible set, the optimal point lies on the boundary of the feasible set for the initial point within the feasible set (Fig. \ref{xeqmout_x0in}) as well as for the initial point lying outside the feasible set (Fig. \ref{xeqmout_x0out}).
 \\
\subsubsection{\textbf{ Convex Objective Function with Inequality Constraints}}
Consider a constrained optimization problem from \cite{hauswirth2024optimization} with a convex objective function with inequality constraints.
\begin{equation}\label{example2}
    \begin{split}
  &\mathrm{min}\; \mathrm{f}(\mathrm{x})\\
  &\mathrm{s.t.}\; \mathrm{x}\in \mathrm{X}\\
  &\mathrm{g}(\mathrm{x})\leq 0
    \end{split}
\end{equation}
where $\mathrm{f(x)=-0.5x_1+0.25x_2}$ with the constraints being $\mathrm{x_2 \geq 0}$ and $\mathrm{x_1 \leq x_2}$. Since the $\mathrm{f(x)}$ and $\mathrm{g(x)}$ are convex and not strictly convex, a penalty function is incorporated to ensure the convergence with inequality constraints. Generally, for the constraints of the form $\mathrm{g(x) \leq 0}$, the commonly used penalty function is the squared 2-norm of the constraint violation vector i.e. $\frac{\rho }{2}\left \| \mathrm{max\left \{ g(x),0 \right \}} \right \|^2$, where $\rho>0$ represents a scaling parameter \cite{hauswirth2024optimization}.  With the penalty function \cite{allibhoy2023control} the above objective function is modified as $\mathrm{f(x)=0.25\left \| x \right \|^2-0.5x_1+0.25x_2}$. 
\begin{remark}
    The penalty function is necessary to ensure that the target dynamics is exponentially stable. Without the penalty function, hessian $\mathrm{\nabla ^2 f(x)}$ is 0, and the condition in (\ref{target exponential}) cannot be fulfilled. Hence, the penalty term is necessary to ensure $\gamma>0\Rightarrow \mathrm{\nabla^2f(x)-k\nabla^T g(x)\nabla g(x)}>0$.
\end{remark}
For a simple representation of the PDGD dynamics, the above optimization problem can be reduced to a matrix notation as minimizing $\mathrm{f(x)=x^TWX+Fx}$ subject to $\mathrm{Ax}\leq 0$, where $\mathrm{W=\begin{bmatrix}
0.25 &0 \\ 
0 &0.25 
\end{bmatrix}}$, $\mathrm{F=\begin{bmatrix}
-0.5 &0.25 
\end{bmatrix}^T}$ and $\mathrm{A=\begin{bmatrix}
1 &-1 \\ 
0 &-1 
\end{bmatrix}}$.

The Controlled PDGD dynamics for the above problem is written as 
\begin{equation} \label{CPDGD2}
    \begin{split}
 \dot{\mathrm{x}}=-\nabla \mathrm{f(x)}-\nabla=\mathrm{-2Wx-F-A^T\lambda} \\\dot{\lambda}=\mathrm{g(x)=\mathrm{Ax+u}}
    \end{split}
\end{equation}

The invariant manifold for the above problem is $\Psi({x,\lambda})=\lambda+\mathrm{kAx}=0$. After the construction of
the invariant manifold, the convergence of the trajectories of
off-the-manifold dynamics (manifold attractivity) is assured
through the control input. The control input is related to a storage function, which in turn is obtained from the splitting (connection) of the
tangent space in the fiber bundle, defining the given dynamical
system using the tangent kernel (TK) (which is a degenerate two-form resembling a pseudo-Riemannian (PR) metric). The (non-quadratic) candidate Lyapunov function is defined as $\mathrm{S(x,\lambda)=\frac{1}{2}(\lambda+kAx)^2}$. With the condition of exponential convergence $\mathrm{\dot{S}\leq -\alpha S}$, the control law is obtained as
\begin{small}
   \begin{equation} \label{Control2}
\mathrm{u=-Ax+2kAWx+kAF+kAA^T\lambda-\frac{\alpha}{2}(\lambda+kAx)}
\end{equation} 
\end{small}
To guarantee that the optimal solution is evaluated while fulfilling the constraints, the projection operator is used. 
\begin{figure}[ht!]
 \centering
     \includegraphics[ width=\linewidth]{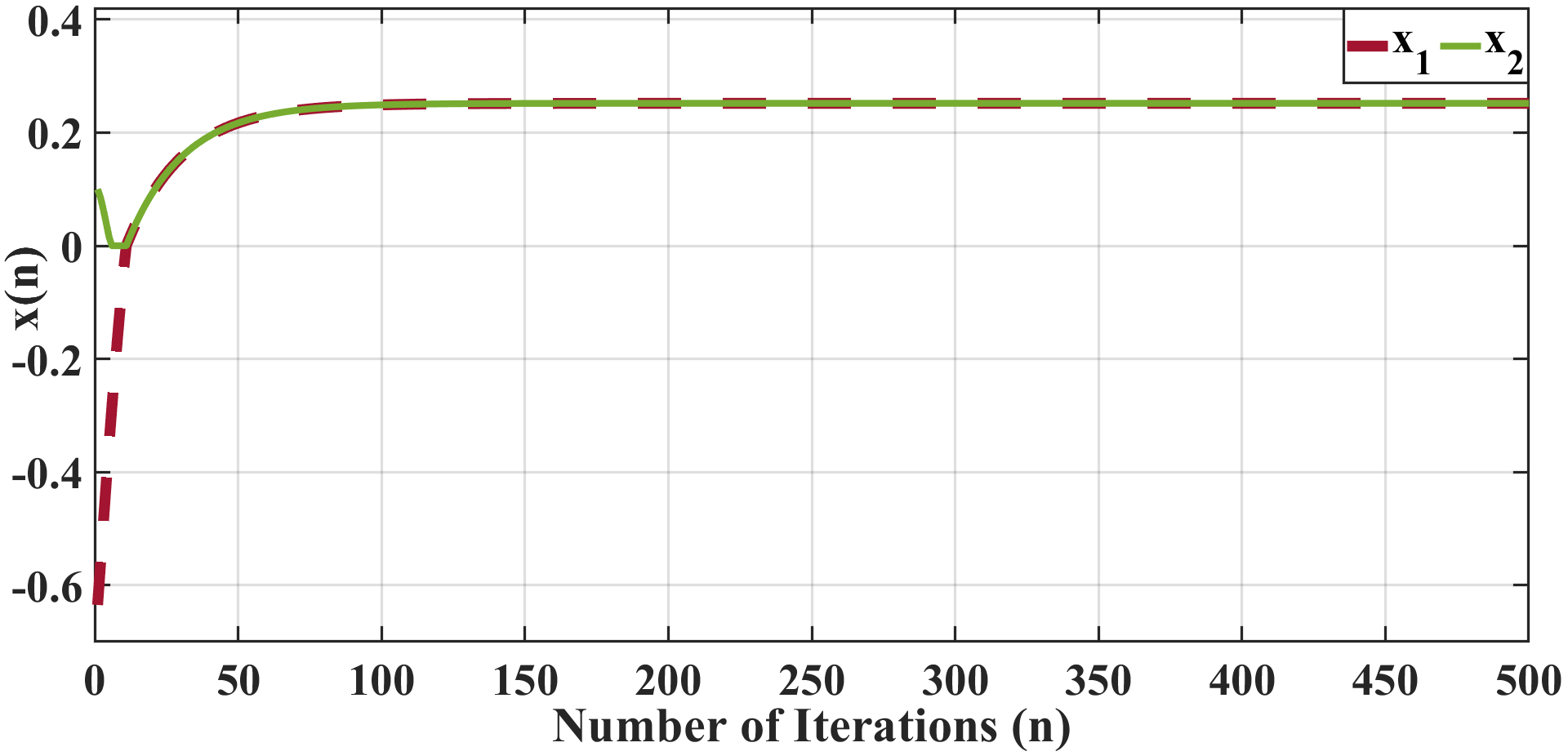}
     \caption{Evolution of the system dynamics in (\ref{CPDGD2}) with the control input (\ref{Control2}) for $\alpha=10$, $k=0.01$ and step size $\beta=0.1$. The y-axis illustrates the values of variables and the x-axis illustrates the number of iterations.}
      \label{cortes_proj}
 \end{figure}
\begin{figure}[ht!] 
 \centering
     \includegraphics[ width=\linewidth]{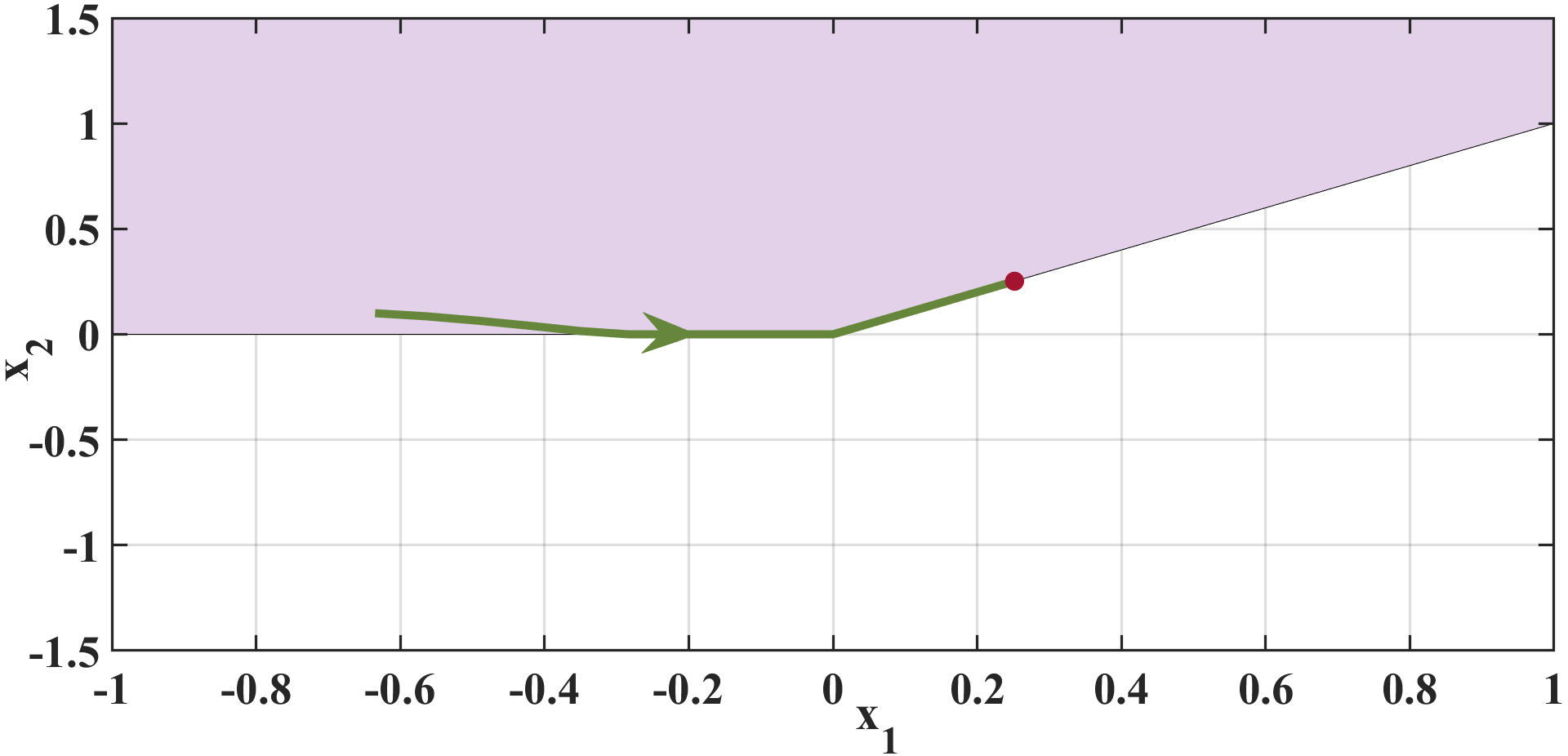}
     \caption{Region of Operation: The purple-shaded region is the feasible set characterized by the constraints $\mathrm{x_1 \leq x_2}$ and $\mathrm{x_2 \geq 0}$. The solution trajectory is shown with the green line and the optimal solution is denoted by the red dot.}
     \label{2Dcortes_proj}
 \end{figure}
The time response of the controlled PDGD dynamics in (\ref{CPDGD2}) with the control law (\ref{Control2}) is illustrated in Fig. \ref{cortes_proj}. With the selected invariant manifold, for the minimum values of $\mathrm{k}$, the minima of $\mathrm{L(x)=f(x)+\lambda^Tg(x)}$ ( within the feasible region) closer to the global minima of the unconstrained objective function $\mathrm{f(x)}$ is obtained. The controlled PDGD improves the performance of the optimizer through the notion of an invariant manifold and its attractivity as compared to the results of various methods mentioned in \cite{allibhoy2023control}. From Fig. \ref{2Dcortes_proj}, it is observed that the trajectory is smooth and it converges to the global minimizer and remains inside the feasible set. Due to the projection on the feasible set, the trajectory remains on the boundary of the feasible region.\\

\subsubsection{\textbf{Nonconvex Objective Function with Nonlinear Constraints}}
The proposed approach is generalized enough to accommodate a variety of constrained optimization problems. In the previous examples, the considered objective functions were convex and the constraints were linear. In this section, a constrained optimization problem with a nonconvex objective function subject to nonlinear constraints is explored.

\begin{figure}[ht]
 \centering
     \includegraphics[ width=\linewidth]{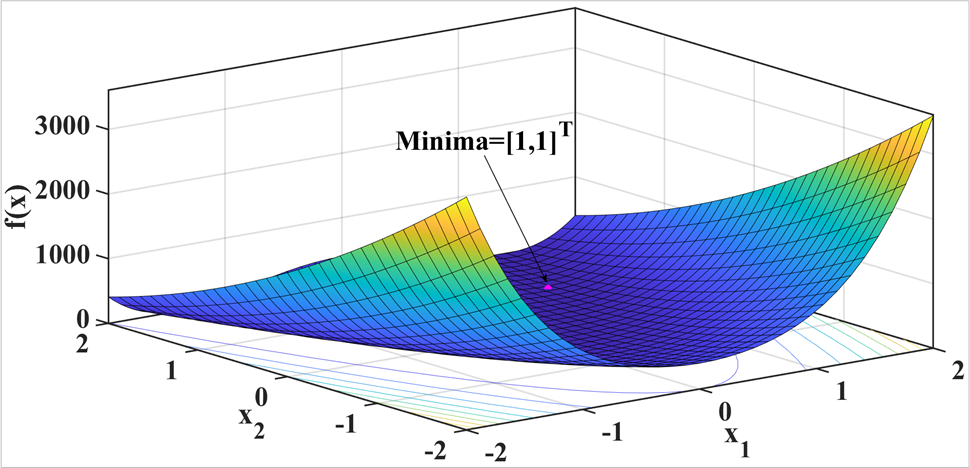}
     \caption{Illustration of Rosenbrock function with minima at $[1,1]^T$}
     \label{mesh1}
 \end{figure}

A Rosenbrock function, which is a nonconvex objective function is considered. Rosenbrock is a benchmark test cost function for global optima-seeking problems. The proposed approach is demonstrated for the Rosenbrock function due to its hard-to-find minima. 
As illustrated in Fig. \ref{mesh1}, the global optimum of the Rosenbrock function exists at $\textcolor{black}{x}^*=\left [ 1,1 \right ]^T$ within a large parabolic-shaped valley, which makes optimizing the Rosenbrock function computationally hard.  
The GD dynamics of the Rosenbrock function is an ill-conditioned system \cite{beck2014introduction} as evident from the higher spectral condition number of the Hessian of Rosenbrock function at its minima. Due to the stiffness of the system, minimizing the Rosenbrock cost function requires more computational cost \cite{bhattacharjee2021closed} and the gradient dynamics often converge slowly to the minima.
\begin{remark}
    The spectral condition number of the Hessian of cost function evaluated at its minima is useful to characterize the performance of the optimization technique and is given by $\kappa =\frac{\lambda _{max}}{\lambda _{min}}$, where $\lambda$ denotes the eigenvalue \cite{gilbert2006linear}. The spectral condition number of the Hessian of Rosenbrock function at the minima $\textcolor{black}{x}^*=\left [ 1,1 \right ]^T$ is calculated as  $\kappa=2508$
\end{remark}
Consider a example from \cite{allibhoy2021anytime}
\begin{equation}\label{example3}
    \begin{split}
&\underset{\mathrm{x}\in\mathbb{R}^2}{\mathrm{minimize}}\;\; (\mathrm{1-x_1})^2+100(\mathrm{x_2-x_1^2})^2\\
&\mathrm{subject\, to}\;\; \mathrm{x_1^2+x_2^2\leq2}
    \end{split}
\end{equation}

The controlled PDGD dynamics for the above optimization problem is of the form
\begin{equation}\label{CPDGD3}
    \begin{split}
&\mathrm{\dot{x}_1=400x_1(x_2-x_1^2)+2-2x_1-2x_1\lambda}\\
&\mathrm{\dot{x}_2=200x_1^2-200x_2-2x_2\lambda}\\
&\mathrm{\dot{\lambda}=x_1^2+x_2^2-2+u}
    \end{split}
\end{equation}
The invariant manifold for the considered problem will be of a form
\begin{equation}
    \Psi(\mathrm{x,\lambda)}=\mathrm{\lambda+k(x_1^2+x_2^2-2)=0}
\end{equation}
The target dynamics defined on the manifold is exponentially stable with the lower values of the parameter $\mathrm{k}$. The attractivity of the manifold is ensured through the control law evaluated through the P\&I approach. The corresponding candidate Lyapunov function is $\mathrm{S(x,\lambda)=\frac{1}{2}(\lambda+k(x_1^2+x_2^2-2))^2}$. Through the exponential convergence condition, the final control law is obtained as
\begin{small}
   \begin{equation}\label{Control3}
   \begin{split}
       &\mathrm{u=-x_1^2-x2^2+2-k(800x_1^2(x_2-x_1^2)+4x_1-4x_1^2-4x_1^2\lambda)}\\-&\mathrm{k(400x_1^2x_2-400x_2^2-4x_2^2\lambda)-\frac{\alpha}{2}(\lambda+k(x_1^2+x_2^2-2))} 
   \end{split}
\end{equation} 
\end{small}
\begin{figure}[ht!] 
 \centering
     \includegraphics[ width=\linewidth]{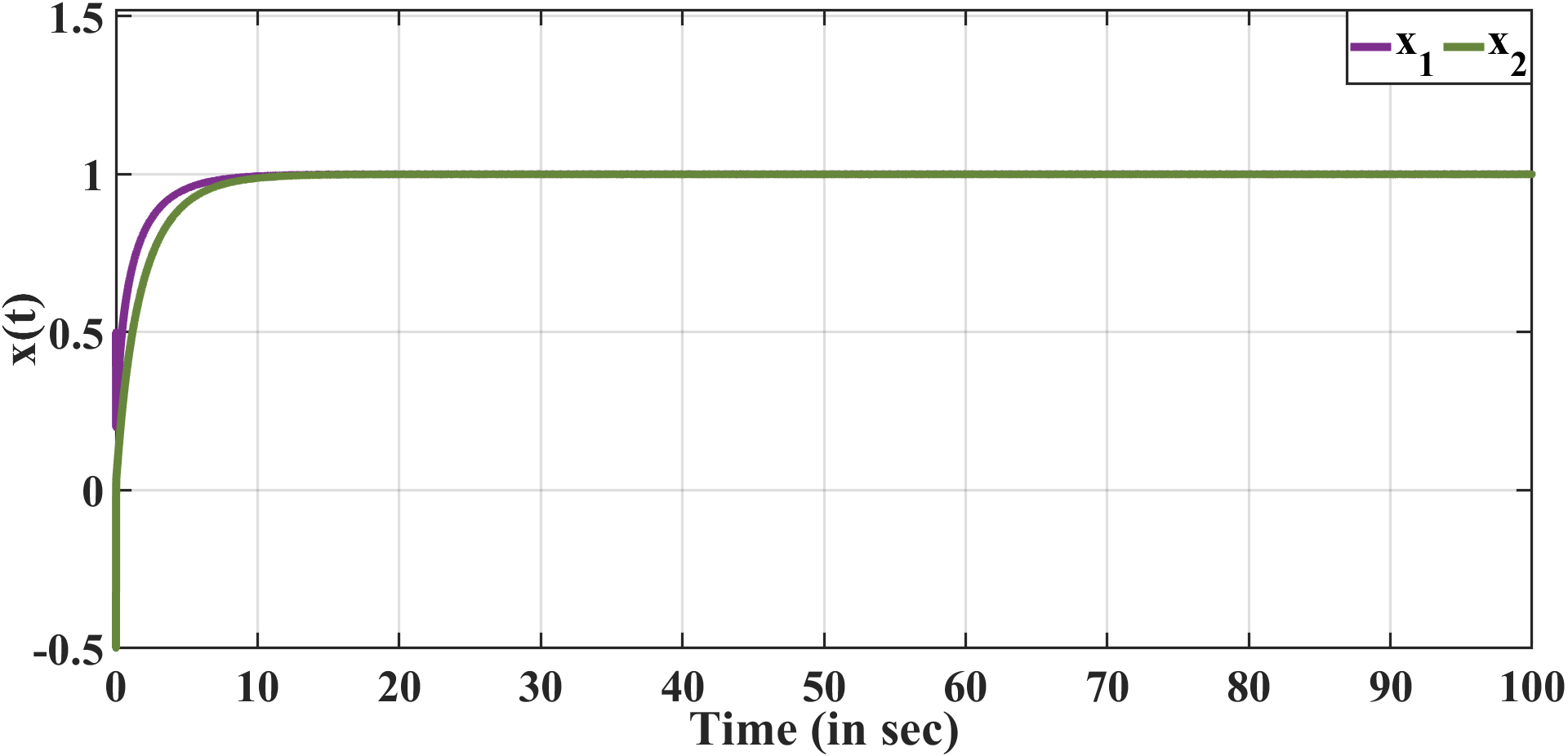}
     \caption{Time evolution of the system in (\ref{CPDGD3}) with the control input (\ref{Control3}) for $\alpha=5$ and $k=0.001$ The y-axis illustrates the values of variables and the x-axis illustrates the time in seconds.}
     \label{RBNonlinear}
 \end{figure}
The above control law ensures that the optimal solution is evaluated while fulfilling the constraints. The time response of the above controlled-PDGD dynamics with the parameter values $\alpha=5$ and $\mathrm{k=0.001}$ are demonstrated in the Fig. \ref{RBNonlinear}. It is observed that irrespective of the nonlinear constraints, the trajectories of controlled PDGD dynamics converge to the global minima while remaining in the feasible region. \\

\subsubsection{\textbf{Example 4: Non-convex Objective Function with Non-convex Constraints}}
The proposed controlled-PDGD framework is applicable for the constrained optimization problem involving the non-convex constraints given that the global minima lies within the feasible region. The proposed framework can be applied for nonconvex constraints $\mathrm{g(x)}$ without any loss of generality. Consider an optimization problem from \cite{simionescu2002new}, with the non-convex Rosenbrock objective function and non-convex constraints. 
\begin{equation}\label{example4}
\begin{split}
&\underset{\mathrm{x}\in\mathbb{R}^2}{\mathrm{minimize}}\;\; (\mathrm{1-x_1})^2+100(\mathrm{x_2-x_1^2})^2\\
&\mathrm{subject\, to}\;\; \mathrm{(x_1-1)^3-x_2+1\leq0}\\
&\;\;\;\;\;\;\;\;\hspace{1cm} \mathrm{x_1+x_2-2\leq0}
\end{split}
\end{equation}
The controlled PDGD dynamics for the above optimization problem is of the form
\begin{equation}\label{CPDGD4}
    \begin{split}
&\mathrm{\dot{x}_1=400x_1(x_2-x_1^2)-2x_1+2-3\lambda_1(x_1-1)^2-\lambda_2}\\
&\mathrm{\dot{x}_2=200x_1^2-200x_2+\lambda_1-\lambda_2}\\
&\mathrm{\dot{\lambda}_1=(x_1-1)^3-x_2+1+u_1} \\
&\mathrm{\dot{\lambda}_2=x_1+x_2-2+u_2}
    \end{split}
\end{equation}
The control law $\mathrm{u_1}$ is evaluated through the systematic procedure of P\&I using the invariant manifold $\mathrm{\Psi(x,\lambda)=\lambda_1+k((x_1-1)^3-x_2-1)=0}$.
\begin{small}
    \begin{equation}\label{Control4I}
        \begin{split}
            &\mathrm{u_1=-(x_1-1)^3+x_2-1-3k(x_1-1)^2(400x_1(x_2-x_1^2))}\\
            &\mathrm{-3k(x_1-1)^2(-2x_1+2-3\lambda_1(x_1-1)^2-\lambda_2)}\\
            &\mathrm{-k(200x_1^2-200x_2+\lambda_1-\lambda_2)-\frac{\alpha}{2}(\lambda_1+k((x_1-1)^3-x_2+1))}
        \end{split}
    \end{equation}
\end{small}
Similarly, the control law $\mathrm{u_2}$ is evaluated using the invariant manifold $\mathrm{\Psi(x,\lambda)=\lambda_2+k(x_1+x_2-2)=0}$.
\begin{small}
    \begin{equation}\label{Control4II}
        \begin{split}
            &\mathrm{u_2=-x_1-x_2+2-k(400x_1(x_2-x_1^2)-2x_1+2-3\lambda_1(x_1-1)^2-\lambda_2)}\\
            &\mathrm{-k(200x_1^2-200x_2+\lambda_1-\lambda_2)-\frac{\alpha}{2}(\lambda_2+k(x_1+x_2-2))}
        \end{split}
    \end{equation}
\end{small}
\begin{figure}[ht] 
 \centering
     \includegraphics[ width=\linewidth]{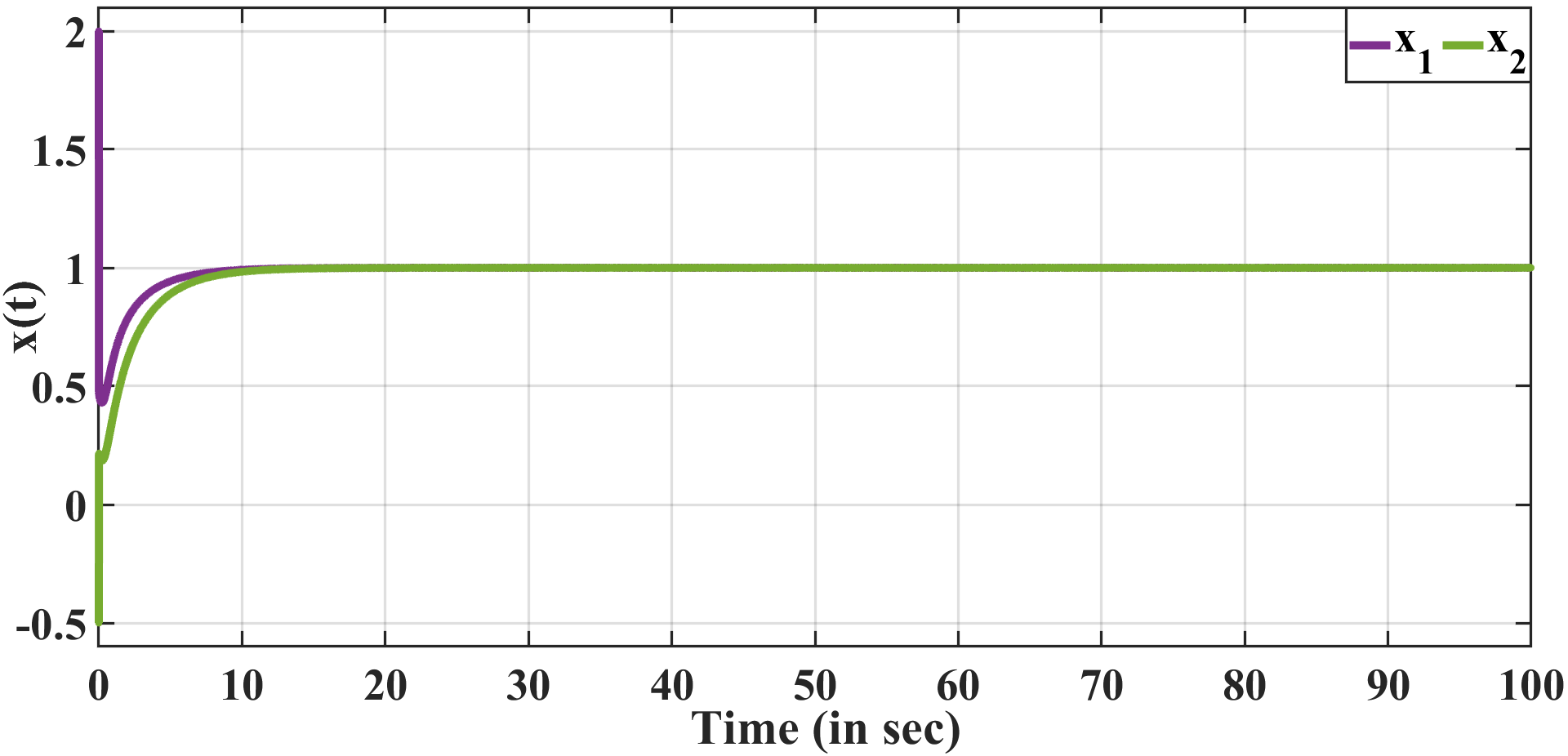}
     \caption{Time evolution of the system in (\ref{CPDGD4}) with the control inputs (\ref{Control4I}) and (\ref{Control4II}) for $\alpha=5$ and $k=0.01$ The y-axis illustrates the values of variables and the x-axis illustrates the time in seconds.}
     \label{RBnonconvex}
 \end{figure}
The time response for the above Controlled PDGD dynamics along with the control laws (\ref{Control4I}), (\ref{Control4II}) and the parameters $\alpha=5$, $k=0.01$ is illustrated in Fig. \ref{RBnonconvex}. The trajectories converge to the global minima while satisfying the constraints.

 The proposed framework is inadequate for the constrained optimization problems with the nonconvex constraints and the global minima lying outside the feasible region, as the projection on the feasible set approach fails in such cases. To alleviate these limitations, a different invariant manifold should be selected such that the requirement of the projection operator is avoided. 
 
 These limitations have motivated the use of a Control Barrier Function (CBF) inspired invariant manifold. The CBF-inspired invariant manifold is constructed with the philosophy that the optimizer should always stay within the feasible region of operation. The CBF prevents the optimizer from leaving the feasible set by pushing the optimizer within the feasible region of the operation. Instead of selecting an invariant manifold of the form $\mathrm{\Psi(x,\lambda)=\lambda+kg(x)=0}$, another modified manifold of the form $\mathrm{\Psi(x,\lambda)=\lambda^2+kg(x)=0}$ should be selected. This modified invariant manifold acts like a CBF and satisfies the conditions (\ref{dual constraints}) and ($\ref{constraints condition}$). 
\section{A Control Perspective in Optimization-based Control} \label{control in CBF}
\subsection{Limitations of the projection-based controlled PDGD dynamics}
The continuous projection defined in (\ref{projection operator}) projects the dynamics onto the feasible (or constraint) set. However, when the constraint function is non-convex, the projection methods become insufficient to ensure the feasibility of the decision variable. In such scenarios, the above-controlled PDGD dynamics should be modified. Hence, the implicit manifold can be designed in such a way that the stabilization of the dynamics is ensured along with the safety objective (i.e. to remain in the feasible region).
\subsection{CBF Inspired Controlled PDGD}
The controlled PDGD is modified for accommodating the constraint optimization problems with the non-convex constraints. The modification is inspired by the notion of the control barrier function.

Taking inspiration from the idea of the control barrier function (CBF) \cite{allibhoy2023control} \cite{ames2016control}, an invariant manifold is defined as 
\begin{equation}\label{CBF_IM}
    \mathrm{\Psi(x,\lambda)=\lambda^2+kg(x)=0}
\end{equation}
This manifold will in turn act as a CBF ensuring that the optimizer (GD algorithm) will not leave the feasible region. When the optimizer enters the zone outside the feasible set, the CBF will push the optimizer inside the feasible set to ensure the safety objectives of the operations. 
As per the Lemma 1, the above manifold is invariant. After the construction of the invariant manifold, its attractivity is assured through the control law evaluated by following a systematic procedure of P\&I. Through the splitting (connection) of the tangent bundle with the PR metric, a storage function $\mathrm{S(x,\lambda)=\frac{1}{2}(\lambda^2+kg(x)^2}$ is obtained. The control law is formulated using the condition of exponential convergence, i.e., $\mathrm{\dot{S}\leq -\alpha S}$
\begin{small}
   \begin{equation}\label{CBF_Control}
\mathrm{u=-g(x)-\frac{1}{2\lambda}\left ( k\nabla g(x)(-\nabla f(x)-\nabla ^T g(x) \lambda)+\frac{\alpha}{2}(\lambda^2+kg(x)) \right )}
\end{equation} 
\end{small}
Hence, with the CBF-inspired invariant manifold the controlled PDGD dynamics is modified as
\begin{small}
  \begin{equation}\label{CBF_PDGD_with_control}
    \begin{split}
\dot{\mathrm{x}}=-\nabla \mathrm{f(x)}-\nabla \mathrm{^Tg(x)}\lambda\\
\dot{\lambda}=\mathrm{-\frac{1}{2\lambda}\left ( k\nabla g(x)(-\nabla f(x)-\nabla ^T g(x) \lambda)+\frac{\alpha}{2}(\lambda^2+kg(x)) \right )}
    \end{split}
\end{equation}  
\end{small}
\subsection{Numerical Example: CBF-Inspired Controlled PDGD}
To understand the applicability of the CBF-inspired invariant manifold, let's revisit the earlier example with the convex objective function and inequality constraints from \cite{hauswirth2024optimization}. 

\begin{equation}\label{examplecbf}
    \begin{split}
  &\mathrm{min}\; \mathrm{x^TWX+Fx}\\
  &\mathrm{s.t.}\; \mathrm{Ax}\leq 0\\
    \end{split}
\end{equation}
where $\mathrm{W=\begin{bmatrix}
0.25 &0 \\ 
0 &0.25 
\end{bmatrix}}$, $\mathrm{F=\begin{bmatrix}
-0.5 &0.25 
\end{bmatrix}^T}$ and $\mathrm{A=\begin{bmatrix}
1 &-1 \\ 
0 &-1 
\end{bmatrix}}$. The corresponding PDGD dynamics along with the control is as follows:
\begin{equation} \label{CPDGDcbf}
    \begin{split}
 \dot{\mathrm{x}}=-\nabla \mathrm{f(x)}-\nabla=\mathrm{-2Wx-F-A^T\lambda} \\\dot{\lambda}=\mathrm{g(x)=\mathrm{Ax+u}}
    \end{split}
\end{equation}

With reference to (\ref{CBF_IM}), the invariant manifold is written as $\mathrm{\Psi(x,
\lambda)=\lambda^2+kAx=0}$. With this manifold and following the step-wise procedure of the P\&I approach, the control law is evaluated as 
\begin{small}
    \begin{equation}\label{CBF_control_law}
        \mathrm{u=-Ax-\frac{1}{2\lambda}\left ( -2kAWx-kAF-kAA^T\lambda+\frac{\alpha}{2}(\lambda^2+kAx) \right )}
    \end{equation}
\end{small}
\begin{figure}[ht!] \label{CBF_cortes_x}
 \centering
     \includegraphics[ width=\linewidth]{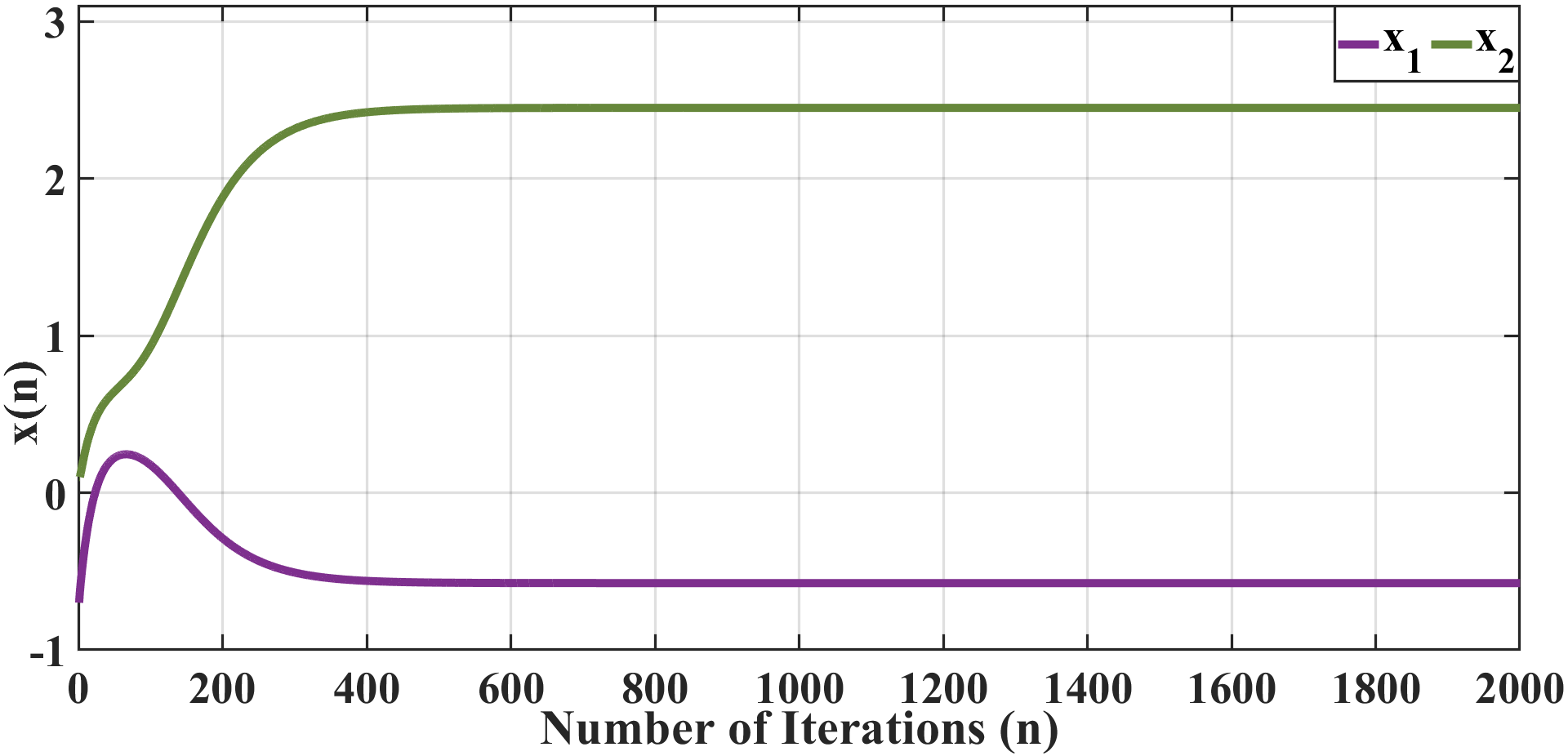}
     \caption{Evolution of the system dynamics in (\ref{CPDGDcbf}) with the control input (\ref{CBF_control_law}) for $\alpha=10$, $k=0.2$ and step size $\beta=0.1$. The y-axis illustrates the values of variables and the x-axis illustrates the number of iterations.}
 \end{figure}
\begin{figure}[ht!] \label{CBF_cortes_u}
 \centering
     \includegraphics[ width=\linewidth]{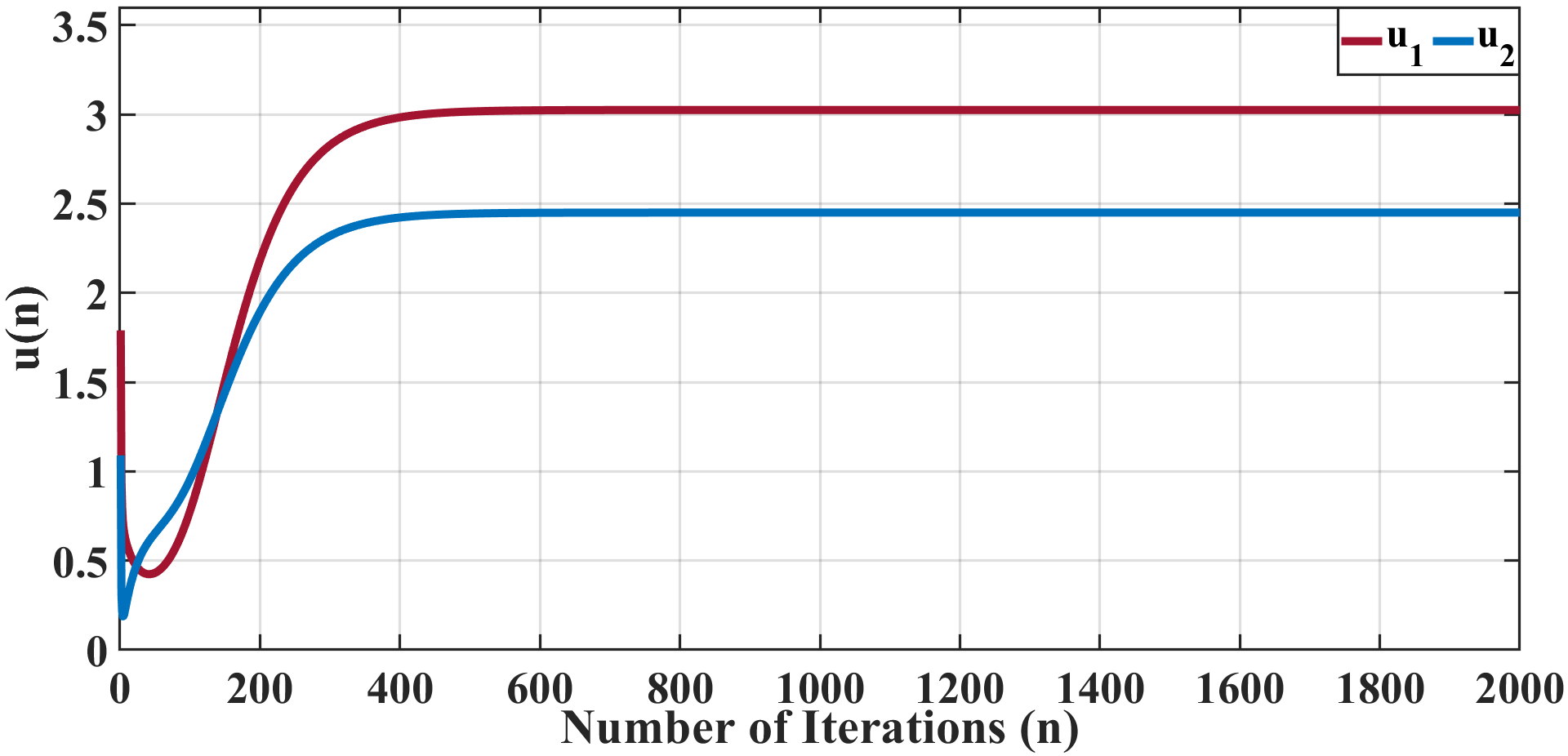}
     \caption{Evolution of the control input (\ref{CBF_control_law}) for $\alpha=10$, $k=0.2$ and step size $\beta=0.1$. The y-axis illustrates the values of control inputs and the x-axis illustrates the number of iterations.}
 \end{figure}
 \begin{figure}[ht!] \label{2Dcortes_CBF}
 \centering
     \includegraphics[ width=\linewidth]{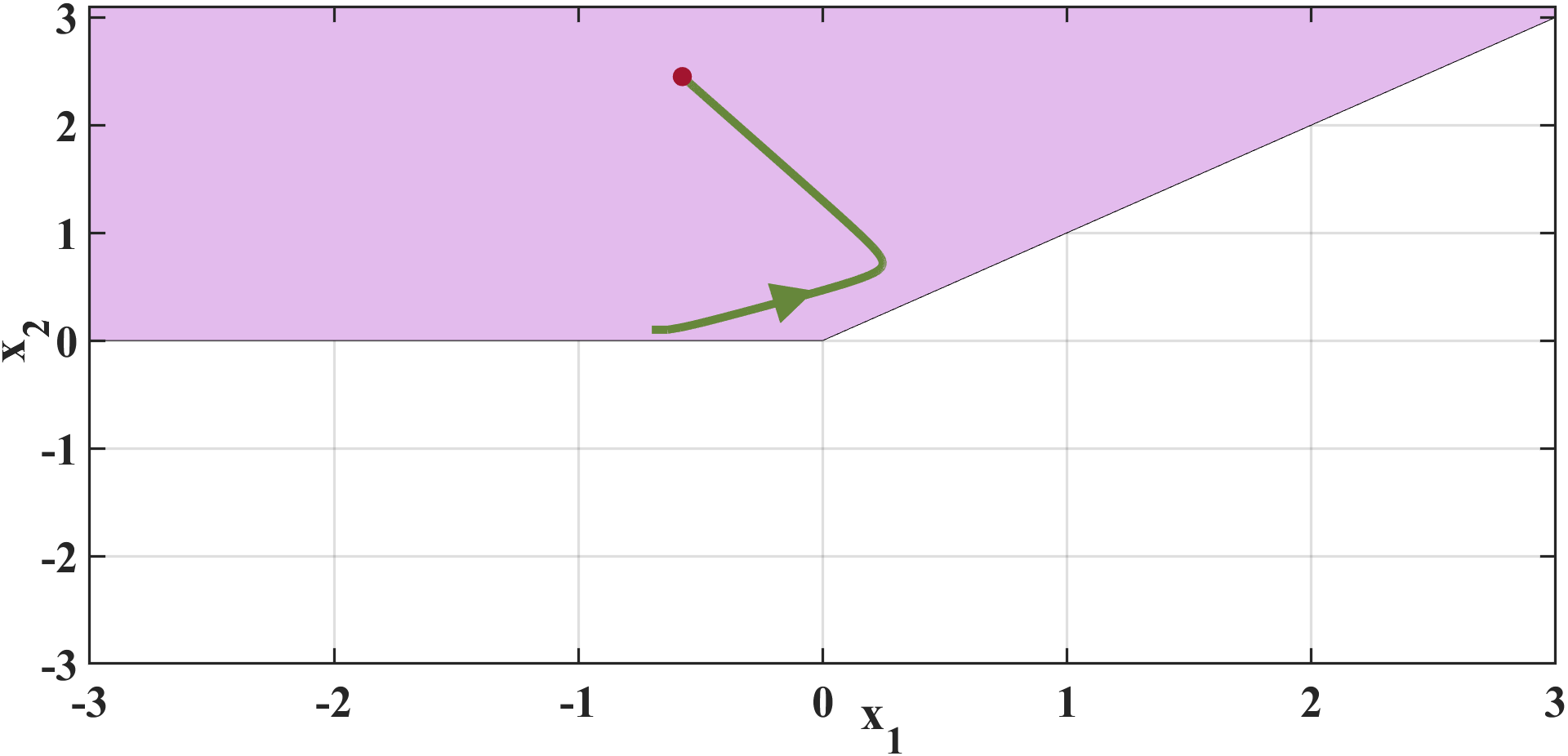}
     \caption{Region of Operation: The purple-shaded region is the feasible set characterized by the constraints $\mathrm{x_1 \leq x_2}$ and $\mathrm{x_2 \geq 0}$. The solution trajectory is shown with the green line and the optimal solution is denoted by the red dot.}
 \end{figure}
The time response of the system dynamics (\ref{CPDGDcbf}) with the control input (\ref{CBF_control_law} are illustrated in Fig.\ref{CBF_cortes_x} and Fig.\ref{CBF_cortes_u} respectively. The region of operation with the feasible set is represented in Fig.\ref{2Dcortes_CBF}. It is evident from Fig.\ref{2Dcortes_CBF} that the invariant manifold acts as a CBF and pushes the optimizer inside the feasible region whenever it tends to leave the feasible set. The evolution of the control inputs is represented in Fig. \ref{CBF_cortes_u}, which pushes the optimizer inside the feasible region and hence the GD dynamics converge to the equilibrium point inside the feasible region. The optimal point can be further improved with various values of the parameter $k$ and the choice of an invariant manifold.
\begin{remark}
     The invariant manifold of the quadratic nature is selected in the above formulation. However, the manifolds of logarithmic nature $\mathrm{\Psi(x,\lambda)=log(\lambda)+kg(x)=0}$ can also be selected. With this manifold, the target dynamics is written as $\mathrm{\dot{x}=-\nabla f(x)-e^{-kg(x)}\nabla g(x)}$. The effect of $\mathrm{\nabla g(x)}$ remains undiminished and hence $\mathrm{x}^*$ farther from the global minima is achieved. The optima achieved with the manifold $\mathrm{\Psi(x,\lambda)=\lambda^2+kg(x)=0}$ is better than the logarithmic manifold. 
\end{remark}

 The above-defined formulation is useful in optimization-based control applications where safety (stabilization within the feasible region) is the priority. For example, the integral windup in the PID controller. The optimization-based control applications often involve a quadratic programming problem formulated with the Control Lyapunov Function (CLF) and Control Barrier Function (CBF). Finding the CLF as well as CBF is a tedious task. Hence, in such scenarios, the proposed CBF-inspired Controlled PDGD approach is feasible, since it avoids the formulation of a Quadratic Programming problem saving the computation cost. The CBF-inspired invariant manifold provides the non-quadratic storage function through the systematic procedure of P\&I. This storage function is a `\textit{`Manifold Control Lyapunov Function"} since it is constructed to achieve the attractivity of the manifold. \textcolor{blue}{The proposed formulation incorporates the functionalities of the CLF and CBF without the necessity of solving a Quadratic Programming problem, which formulates another significant contribution of the paper.}

This approach is beneficial in constrained optimization problems such as optimal power flow, where the constraints are non-convex and various relaxation techniques are used to approximate the Lagrangian function to ensure the convergence of the solution. \textcolor{blue}{In such applications, the CBF-inspired controlled PDGD approach proves beneficial as it achieves the optimal solution while satisfying the constraints without any necessity for relaxation techniques.}

\section{A Control Perspective in Parameter Estimation} \label{Control in PE}
\subsection{Controlled Gradient Estimator (CGE)}
The development of various system identification and adaptive control techniques commences with the formulation of a Linear Regression Equation (LRE) which is linear in the unknown parameters of the plant or the controller. 
\begin{equation}\label{LRE}
    \mathrm{y(t)=\phi^T(t)\theta+\epsilon(t)}
\end{equation}
where $y(t)\in \mathbb{R}$ denotes the measurable output signal of the system, $\phi(t)\in\mathbb{R}^q$ represents the measurable known regressor vector, $\epsilon(t)$ is a (generic) exponentially decaying signal, and $\theta\in \mathbb{R}^q$ is a constant vector of the unknown parameters of the system. Consider the online estimate of the parameter $\theta$ is denoted as $\hat{\theta}$, and the corresponding estimates of the output are given as $\mathrm{\hat{y}=\phi^T\hat{\theta}}$. With the ideal values and estimates of $\theta$ and$\mathrm{y}$, the output error equation is formulated as 
\begin{equation}
    \mathrm{\tilde{y}=\phi^T(t)\tilde{\theta}}
\end{equation}
From the above formulation, it is intuitive that, when $\mathrm{\hat{\theta}\rightarrow \theta}$, the output error will diminish, i.e., $\mathrm{\tilde{y}\rightarrow 0}$ as $\mathrm{t\rightarrow 0}$. Hence, the objective is to develop an update law for the unknown parameters using the known regressor vector and output signals. The classical choice to accomplish this objective is a Gradient Descent (GE) technique through which the quadratic error function $\frac{1}{2}(\mathrm{y-\phi^T\hat{\theta}})^2$ is minimized. The GD estimator for the unknown parameters of (\ref{LRE}) is of the form
\begin{equation} \label{GDPE}
    \mathrm{\dot{\hat{\theta}}=-\gamma \phi(y-\phi^T\hat{\theta})}
\end{equation}
where $\gamma$ is the adaptation gain. From the above equation and the parameter error, a well-known Parametric Error Equation (PEE) \cite{ortega2020modified} is formulated as
\begin{equation}\label{PEE}
    \mathrm{\dot{\tilde{\theta}}=-\gamma \phi\phi^T\tilde{\theta}}
\end{equation}
which is an LTV system. The Gradient Estimator (GE) is inadequate to ensure the convergence of the parameters when the regressor vectors are not fulfilling the PE condition. 
\begin{defi}
    A bounded vector signal $\phi(t) \in \mathbb{R}^q$ is persistently exciting (PE), if there exists $\mathrm{T>0}$, $\alpha>0$ such that
    \begin{equation*}
        \mathrm{\int_{t}^{t+\tau}\phi(\tau)\phi^T(\tau) d\tau \geq \alpha I_q, \;\;\; \forall t\geq0}
    \end{equation*}
\end{defi}
Also, in GE to ensure the faster convergence of the parameters, a larger step size cannot be selected, since, the larger step size might miss the global optima. 
\subsubsection{\textbf{CGE for two unknown parameters}}
For an easier understanding of the formulation of a Controlled Gradient Estimator (CGE), a system with $\mathrm{q=2}$ is considered. For $\mathrm{q=2}$, the LRE in (\ref{LRE}) is written as
\begin{equation}
    \mathrm{y=\phi_1\theta_1+\phi_2\theta_2}
\end{equation}
Multiplying throughout with $\phi=[\phi_1 \phi_2]^T$ and applying a first order filter $\mathrm{\mathbb{L}(s)=\frac{1}{s+1}}$ yields a matrix form as
\begin{equation}\label{new LRE}
\mathrm{\begin{bmatrix}
Y_1\\ 
Y_2
\end{bmatrix}=\begin{bmatrix}
\Omega_{11}  & \Omega_{12} \\ 
\Omega_{21} & \Omega_{22}
\end{bmatrix}\begin{bmatrix}
\theta_1\\ 
\theta_2
\end{bmatrix}\Rightarrow Y=\Omega \theta }
\end{equation}
where $\mathrm{\begin{bmatrix}
Y_1\\ 
Y_2
\end{bmatrix}=\begin{bmatrix}
\mathbb{L}(s)(\phi_1 y)\\ 
\mathbb{L}(s)(\phi_2 y)
\end{bmatrix}}$ and $\mathrm{\begin{bmatrix}
\Omega_{11} & \Omega_{12}\\ 
\Omega_{21} & \Omega_{22}
\end{bmatrix}=\begin{bmatrix}
\mathbb{L}(s)(\phi_1\phi_1) & \mathbb{L}(s)(\phi_1\phi_2) \\ 
\mathbb{L}(s)(\phi_2\phi_1) & \mathbb{L}(s)(\phi_2\phi_2)
\end{bmatrix}}$. With the modified LRE in (\ref{new LRE}), the parameter update law and the PEE can be rewritten as $\mathrm{\dot{\hat{\theta}}=-\gamma (Y-\Omega\hat{\theta})}$ and $ \mathrm{\dot{\tilde{\theta}}=-\gamma \Omega\tilde{\theta}}$ respectively \cite{9121700}.

For $\mathrm{q=2}$, the above PEE can be written as individual parametric error dynamics as follows:
\begin{equation}
    \begin{split}
        &\mathrm{\dot{\tilde{\theta}}_1=-\gamma\Omega_{11}\tilde{\theta}_1-\gamma\Omega_{12}\tilde{\theta}_2}\\
&\mathrm{\dot{\tilde{\theta}}_2=-\gamma\Omega_{21}\tilde{\theta}_1-\gamma\Omega_{22}\tilde{\theta}_2}\\
    \end{split}
\end{equation}
\textcolor{blue}{To improve the transient response of the parameter estimator a control perspective is explored in the parameter estimation problem. Through the notion of an invariant manifold and its attractivity achieved through the P\&I approach, the Controlled Gradient Estimator (CGE) is formulated.} The CGE dynamics of PEE is of the form
\begin{equation}
    \begin{split}
          &\mathrm{\dot{\tilde{\theta}}_1=-\gamma\Omega_{11}\tilde{\theta}_1-\gamma\Omega_{12}\tilde{\theta}_2}\\
&\mathrm{\dot{\tilde{\theta}}_2=-\gamma\Omega_{21}\tilde{\theta}_1-\gamma\Omega_{22}\tilde{\theta}_2+u}\\
    \end{split}
\end{equation}
The objective is to improve the transient response of the estimator to achieve a faster convergence of the parameter estimates to the actual parameters of the system. This objective is accomplished through the manifold defined as 
\begin{equation}
    \mathrm{\Psi(\tilde{\theta})=\tilde{\theta}_1-\beta\tilde{\theta}_2=0 }
\end{equation}

\begin{remark}
    The target dynamics $\mathrm{\dot{\tilde{\theta}}_1=-\gamma\Omega_{11}\tilde{\theta}_1-\gamma\beta\Omega_{12}\tilde{\theta}_1}$ is defined on the manifold $\mathrm{\Psi(\tilde{\theta})}$. With the lower values of $\beta$, the target dynamics defined on the manifold exhibit exponential behavior. 
\end{remark}
According to Lemma 1, the given manifold is invariant. The attractivity of the trajectories of the off-the-manifold dynamics to the invariant manifold is achieved through the control law evaluated by following the systematic procedure of P\&I. Through the notion of the splitting of the tangent bundle with the PR metric the candidate Lyapunov function is formulated as $\mathrm{\dot{S}=\frac{1}{2}(\tilde{\theta}_1-\beta\tilde{\theta}_2)^2}$. Through the condition $\mathrm{\dot{S}\leq0}$, the control law is evaluated as:
\begin{small}
    \begin{equation}
        \mathrm{u=\gamma \Omega_{21}\tilde{\theta}_1+\gamma \Omega_{22}\tilde{\theta}_2-\beta\gamma \Omega_{11}\tilde{\theta}_1-\beta\gamma \Omega_{12}\tilde{\theta}_2}
    \end{equation}
\end{small}
To obtain the parameter update law, the above equation can be rewritten by expanding the PEE equations.
\begin{small}
    \begin{equation}
        \mathrm{u=\gamma(Y_2-\Omega_{21}\hat{\theta}_1-\Omega_{22}\hat{\theta}_2)-\gamma\beta(Y_1-\Omega_{11}\hat{\theta}_1-\Omega_{12}\hat{\theta}_2) }
    \end{equation}
\end{small}
where $\mathrm{Y_1=\Omega_{11}\theta_1+\Omega_{12}\theta_2}$ and $\mathrm{Y_2=\Omega_{21}\theta_1+\Omega_{22}\theta_2}$.  
With the above control law, the final Controlled Gradient Estimator model is obtained as:
\begin{small}
    \begin{equation}
    \begin{split}
        &\dot{\hat{\theta}}_1=\gamma(Y_1-\Omega_{11}\hat{\theta}_1-\Omega_{12}\hat{\theta}_2)\\
&\dot{\hat{\theta}}_2=2\gamma(Y_2-\Omega_{21}\hat{\theta}_1-\Omega_{22}\hat{\theta}_2)-\gamma\beta(Y_1-\Omega_{11}\hat{\theta}_1-\Omega_{12}\hat{\theta}_2)
    \end{split}
\end{equation}
\end{small}
The above CGE model can be written in the form of a matrix as
\begin{equation} \label{CGE_Estimator_q2}
    \Rightarrow \dot{\hat{\theta}}=-\gamma\mathcal{P}\left [ Y-\Omega\hat{\theta} \right ]
\end{equation}
where $\mathcal{P}=\begin{bmatrix}
1 &0 \\ 
-\beta & 2
\end{bmatrix}$, $\mathrm{Y=\begin{bmatrix}
Y_1\\ 
Y_2
\end{bmatrix}}$, $\Omega=\begin{bmatrix}
\Omega_{11} & \Omega_{12} \\ 
\Omega_{21} & \Omega_{22}
\end{bmatrix}$ and  $\mathrm{\hat{\theta}=\begin{bmatrix}
\hat{\theta}_1\\ 
\hat{\theta}_2
\end{bmatrix}}$.

To avoid the loss of generality, the PEE estimator is also rewritten by following the notion of error defined as $\tilde{\theta}=\theta-\tilde{\theta}$ and the corresponding PEE gradient $\dot{\tilde{\theta}}=-\dot{\hat{\theta}}$.
\begin{small}
    \begin{equation}
        \begin{split}
    \dot{\tilde{\theta}}_1=-\gamma\Omega_{11}\tilde{\theta}_1-\gamma\Omega_{12}\tilde{\theta}_2\\
\dot{\tilde{\theta}}_2=-2\gamma\Omega_{21}\tilde{\theta}_1-2\gamma\Omega_{22}\tilde{\theta}_2+\gamma\beta\Omega_{11}\tilde{\theta}_1+\gamma\beta\Omega_{12}\tilde{\theta}_2
        \end{split}
    \end{equation}
\end{small}
The above CGE for PEE can be written in the matrix form as:
\begin{equation}\label{CGE_PEE}
    \underbrace{\begin{bmatrix}
\dot{\tilde{\theta}}_1\\ 
\dot{\tilde{\theta}}_1
\end{bmatrix}}_{\dot{\tilde{\theta}}}=-\gamma\underbrace{\begin{bmatrix}
1 &0 \\ 
-\beta & 2
\end{bmatrix}}_{\mathcal{P}}\underbrace{\begin{bmatrix}
\Omega_{11} &\Omega_{12}  \\ 
 \Omega_{21} & \Omega_{22} 
\end{bmatrix}}_{\Omega}\underbrace{\begin{bmatrix}
\tilde{\theta}_1\\ 
\tilde{\theta}_2
\end{bmatrix}}_{\theta}
\end{equation}
\begin{equation}\label{CGE_PEE_q2}
    \Rightarrow \dot{\tilde{\theta}}=-\gamma\mathcal{P}\Omega\tilde{\theta}
\end{equation}
From the above PEE (\ref{CGE_PEE}), it is evident that through the control perspective, the scaling of the classical gradient estimator model is obtained with the matrix $\begin{bmatrix}
1 &0 \\ 
-\beta & 2
\end{bmatrix}$. In the PEE model (\ref{CGE_PEE}), for very negligible values of the parameter $\beta$, the PEE model is approximated as
\begin{equation}
    \begin{bmatrix}
\dot{\tilde{\theta}}_1\\ 
\dot{\tilde{\theta}}_1
\end{bmatrix}=\begin{bmatrix}
-\Omega_{11} &-\Omega_{12}  \\ 
 -2\Omega_{21} & -2\Omega_{22} 
\end{bmatrix}\begin{bmatrix}
\tilde{\theta}_1\\ 
\tilde{\theta}_2
\end{bmatrix}
\end{equation}
The convergence properties of the PEE of the second parameter $\tilde{\theta}_2$ are improved, which further improves the convergence characteristics of the PEE of the first parameter $\tilde{\theta}_1$ through the definition of the manifold. 

The transient response of the parameter estimator is enhanced through the control perspective. The effect of the scaling matrix is reflected in the eigenvalues of the model. With the control perspective, the classical GE is modified such that the minimum eigenvalues of the CGE model increase which in turn improves the transient response of the estimator. 

\subsubsection{\textbf{CGE for three unknown parameters}}
The proposed CGE is generalized enough to be extended to the systems with a higher number of parameters. To understand the generalization, the CGE for the three parameters is demonstrated further.
For $\mathrm{q=3}$, the output equation is defined as
\begin{equation}
    \mathrm{y=\phi_1\theta_1+\phi_2\theta_2+\phi_3\theta_3}
\end{equation}
The CGE model for the PEE of the above system is written as 
\begin{small}
    \begin{equation}
        \begin{split}
            &\dot{\tilde{\theta}}_1=-\gamma\Omega_{11}\tilde{\theta}_1-\gamma\Omega_{12}\tilde{\theta}_2-\gamma\Omega_{13}\tilde{\theta}_3\\
&\dot{\tilde{\theta}}_2=-\gamma\Omega_{21}\tilde{\theta}_1-\gamma\Omega_{22}\tilde{\theta}_2-\gamma\Omega_{23}\tilde{\theta}_3+\mathrm{u_1}\\
&\dot{\tilde{\theta}}_3=-\gamma\Omega_{31}\tilde{\theta}_1-\gamma\Omega_{32}\tilde{\theta}_2-\gamma\Omega_{33}\tilde{\theta}_3+\mathrm{u_2}\\
        \end{split}
    \end{equation}
\end{small}
For calculating $\mathrm{u_1}$ and $\mathrm{u_2}$ independently, two distinct invariant manifolds $\Psi_1(\tilde{\theta})=\tilde{\theta}_2-\beta\tilde{\theta}_1$ and $\Psi_2(\tilde{\theta})=\tilde{\theta}_3-\beta\tilde{\theta}_1$ are considered. 
Through these manifolds, the CGE model for three parameters is derived as 
\begin{small}
    \begin{equation}
        \begin{split}
            \dot{\hat{\theta}}_1=\gamma \epsilon_1\\
\dot{\hat{\theta}}_2=2\gamma\epsilon_2-\gamma\beta\epsilon_1\\
\dot{\hat{\theta}}_3=2\gamma \epsilon_3-\gamma\beta\epsilon_1
        \end{split}
    \end{equation}
\end{small}
where $\epsilon_1=Y_1-\Omega_{11}\hat{\theta}_1-\Omega_{12}\hat{\theta}_2-\Omega_{13}\hat{\theta}_3$, $\epsilon_2=Y_2-\Omega_{21}\hat{\theta}_1-\Omega_{22}\hat{\theta}_2-\Omega_{23}\hat{\theta}_3$, and $\epsilon_3=Y_3-\Omega_{31}\hat{\theta}_1-\Omega_{32}\hat{\theta}_2-\Omega_{33}\hat{\theta}_3$
In the matrix form the above equation can be written as
\begin{equation} \label{CGE_para_3}
    \dot{\hat{\theta}}=-\gamma \mathcal{P}\left [ Y-\Omega\hat{\theta} \right ]
\end{equation}
where $\mathcal{P}=\begin{bmatrix}
1 &0  &0 \\ 
-\beta &2  &0 \\ 
-\beta & 0 &2 
\end{bmatrix}$.
The corresponding CGE model for PEE is
\begin{small}
    \begin{equation}\label{cge_pee_3}
\underbrace{\begin{bmatrix}
\dot{\tilde{\theta}}_1\\ 
\dot{\tilde{\theta}}_2\\ 
\dot{\tilde{\theta}}_3
\end{bmatrix}}_{\dot{\tilde{\theta}}}=-\gamma\underbrace{\begin{bmatrix}
1 &0  &0 \\ 
-\beta &2  &0 \\ 
-\beta & 0 & 2
\end{bmatrix}}_{\mathcal{P}}\underbrace{\begin{bmatrix}
\Omega_{11} & \Omega_{12} &\Omega_{13} \\ 
 \Omega_{21}&\Omega_{22}  & \Omega_{23}\\ 
 \Omega_{31}& \Omega_{32} & \Omega_{33}
\end{bmatrix}}_{\Omega}\underbrace{\begin{bmatrix}
\tilde{\theta}_1\\ 
\tilde{\theta}_2\\ \tilde{\theta}_3
\end{bmatrix}}_{\tilde{\theta}}
    \end{equation}
\end{small}
\begin{remark}
    The proposed CGE approach is a simplified representation of the Modified Gradient Estimator (MGE) of \cite{nayyer2022passivityPE}. The proposed CGE approach is an extension of the work in \cite{nayyer2022passivityPE} for the system with `n' unknown parameters. Instead of adding a single control input, `n-1' control inputs are added to design a generalized framework of CGE.  
\end{remark}

\subsubsection{\textbf{CGE for `n' unknown parameters}}
\textcolor{blue}{For a system with `$\mathrm{q=n}$' unknown parameters, $\mathrm{n-1}$ control inputs are added to the gradient estimators. These control inputs are evaluated distinctly through the distinct manifolds for each control of the form $\tilde{\theta}_i-\beta\tilde{\theta}_1$ where $i=2,3,\cdots, n$.}
The Controlled Gradient Estimator for `n' parameters and for PEE are written as
    \begin{equation}
    \dot{\hat{\theta}}=-\gamma \mathcal{P}\left [ Y-\Omega\hat{\theta} \right ]
\end{equation}
and 
\begin{equation}
    \dot{\tilde{\theta}}=-\gamma \mathcal{P}\Omega \tilde{\theta}
\end{equation}
respectively, where, the scaling matrix evaluated through the control perspective for $\mathrm{n}$ parameters is of form
\begin{equation}
    \mathcal{P}=\begin{bmatrix}
1 &0  &0  &\cdots   &0 \\ 
-\beta &2  &0   &  \cdots & 0\\ 
-\beta & 0 &2  & \cdots  &0 \\ 
\vdots &  \vdots & \vdots  &\ddots& \vdots\\ 
-\beta & 0 &0  & \cdots  & 2
\end{bmatrix}
\end{equation}
\subsection{Numerical Example: Controlled Gradient Estimator}
\subsubsection{\textbf{Two Parameter CGE}}
Consider a system with two parameters $\theta=col(2,-2)$ and a non-PE regressor 
\begin{equation}
    \phi=\begin{bmatrix}
 1& \mathrm{\frac{sin(t)+cos(t)}{(1+t)^{0.5}}-\frac{sin(t)}{2(1+t)^{1.5}}}
\end{bmatrix}^\mathrm{T}
\end{equation}

\begin{figure}[ht]
 \centering
     \includegraphics[ width=\linewidth]{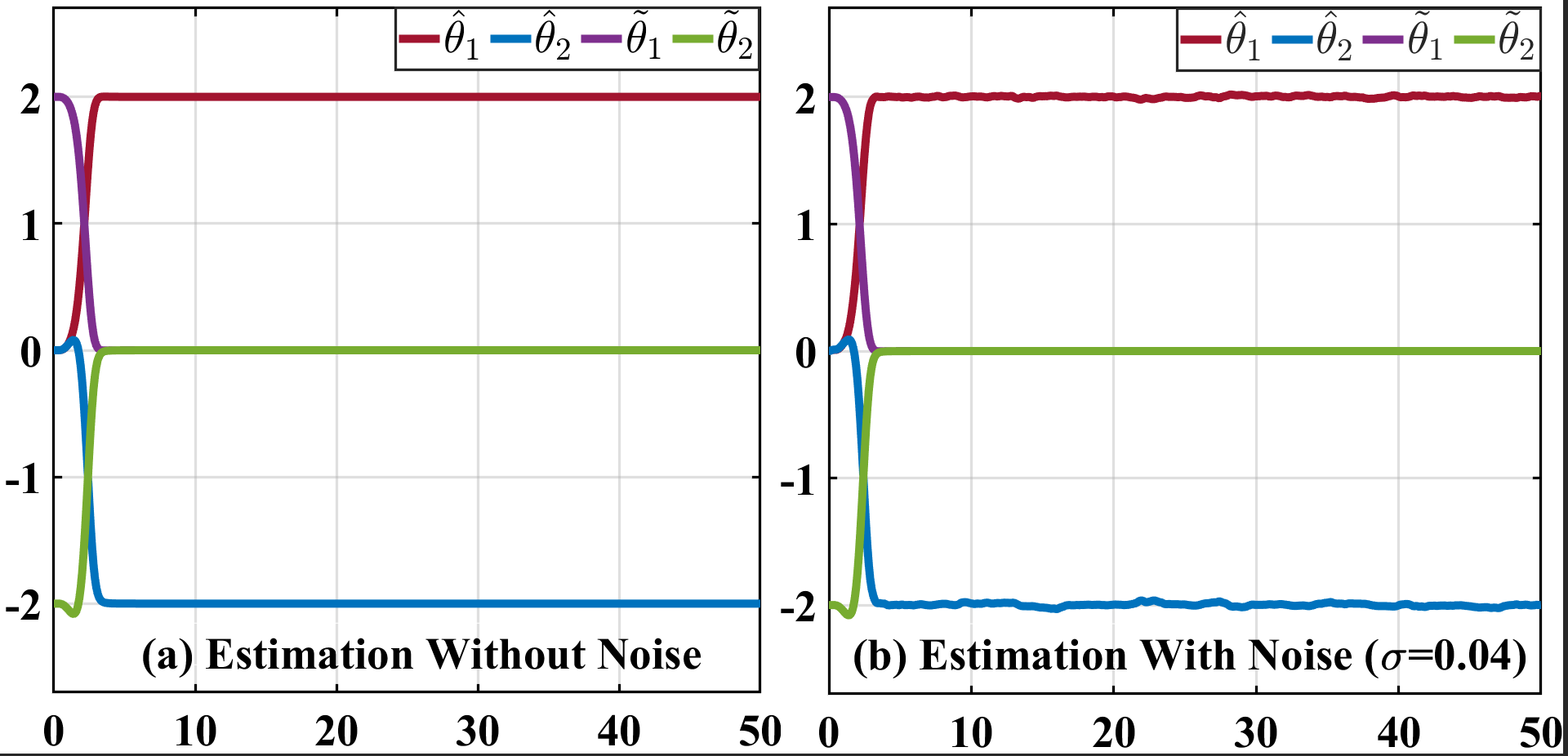}
     \caption{Time evolution of the CGE and the corresponding PEE for $\gamma=10$,  $\beta=0.75$ and $\alpha=10$. Fig. (a) Estimation without noise Fig. (b) Estimation with the measurement noise. The y-axis illustrates the values of variables and the x-axis illustrates the time in seconds.}
     \label{2parape}
 \end{figure}\
 Referring to (\ref{CGE_Estimator_q2}), the Controlled Gradient Estimator for the parameter estimator is formulated. The evolution of CGE $(\hat{\theta}_1(t), \hat{\theta}_2(t)$ with the learning rate $\gamma=10$ and the parameters $\beta=0.75$ and $\alpha=10$ is demonstrated in Fig. \ref{2parape} (a) without any noise in the measurements. The corresponding PEE is formulated with reference to (\ref{CGE_PEE_q2}), and its time response $(\tilde{\theta}_1(t), \tilde{\theta}_2(t))$ has been demonstrated in Fig. \ref{2parape} (a).

 The response of the estimator with the measurement noise (Gaussian noise) 0 mean and 0.04 standard deviation (SD) is illustrated in Fig. \ref{2parape} (b). It is observed that the CGE performs well even in the noisy case due to the notion of the invariant manifold and its attractivity. With the parameter $\alpha$, the convergence rate of the parameter estimates can be adjusted by varying the rate of exponential convergence of the off-the-manifold trajectories. 
\subsubsection{\textbf{Three Parameter CGE}}
 
 \begin{figure}[ht] 
 \centering
     \includegraphics[ width=\linewidth]{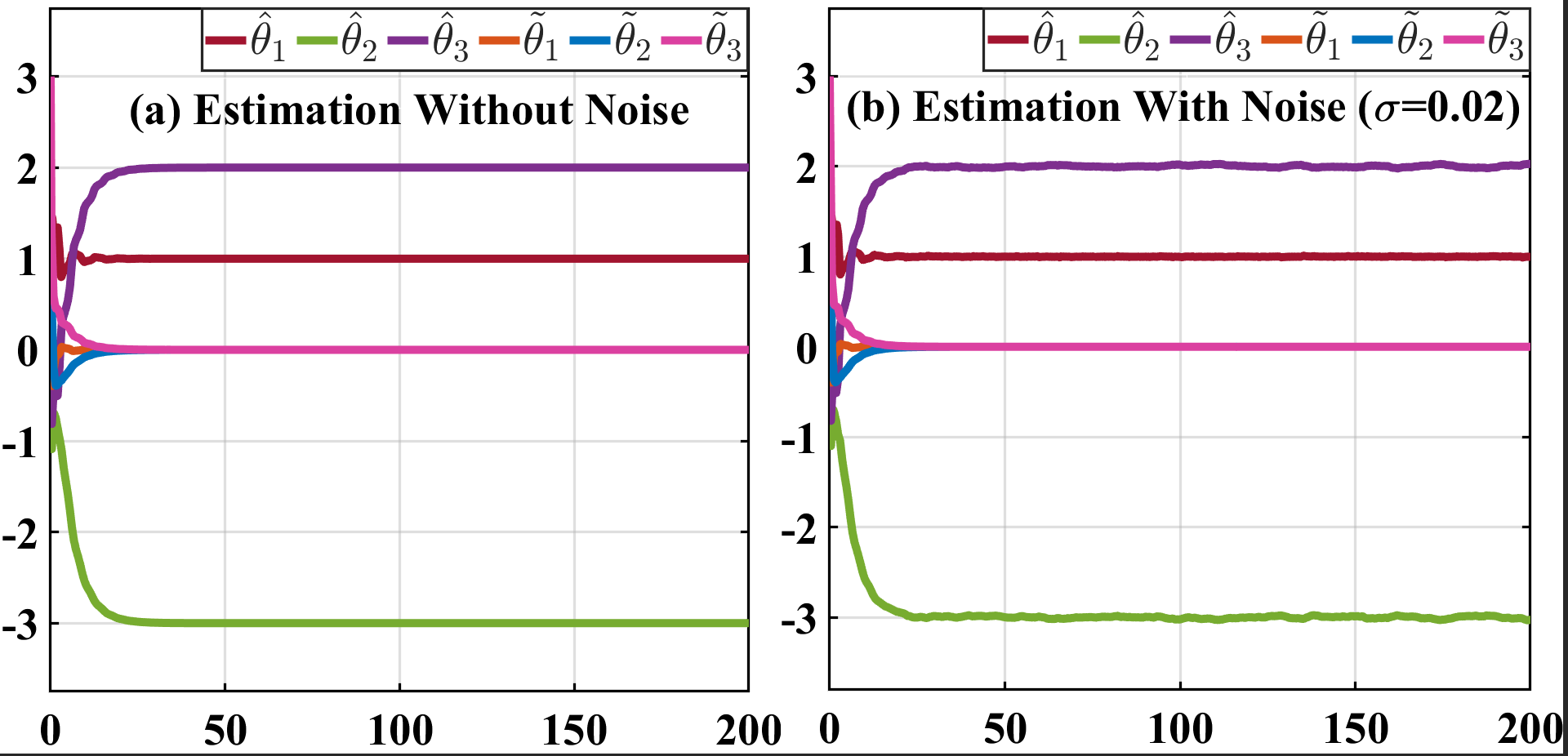}
     \caption{Time evolution of the CGE and the corresponding PEE for $\gamma=100$, and  $\beta=0.95$ and $\alpha=10$. Fig. (a) Estimation without noise Fig. (b) Estimation with the measurement noise. The y-axis illustrates the values of variables and the x-axis illustrates the time in seconds.}
     \label{3parape}
 \end{figure}
The approach has been demonstrated for the system with three parameters $\theta=col(1,2,3)$ with the non PE regressor 
\begin{equation}
    \phi=\begin{bmatrix}
 1&\mathrm{cos(t)}& \mathrm{\frac{sin(t)+cos(t)}{(1+t)^{0.5}}-\frac{sin(t)}{2(1+t)^{1.5}}}
\end{bmatrix}^\mathrm{T}
\end{equation}
For the above example, the parameter update law is formulated with the CGE proposed in (\ref{CGE_para_3}), and the corresponding PEE is formulated referring to (\ref{cge_pee_3}).

 The time response of the parameter estimates of the CGE $(\hat{\theta}_1(t), \hat{\theta}_2(t), \hat{\theta}_3(t))$ and the corresponding PEE $(\tilde{\theta}_1(t), \tilde{\theta}_2(t), \tilde{\theta}_3(t))$ with $\gamma=100$, $\beta=0.95$ and $\alpha=10$ are demonstrated in Fig. \ref{3parape} (a). The performance of the estimator for the noisy case with the Gaussian noise of mean 0 and standard deviation of 0.02 is demonstrated in Fig. \ref{3parape} (b). Due to the notion of the attractivity of the manifold, faster convergence is achieved even in the noisy case. 
\section{Conclusion}
\label{conclusion}
The paper has unified the major themes encompassing an optimization framework such as constrained optimization, optimization-based control, and parameter estimation under one roof through the control perspective.
Rather than making arbitrary modifications to the algorithms through augmentation, incorporating regularization, or adding multiple filters, we are systematically modifying the algorithm by incorporating geometric information through control. Since the problems are represented as dynamical systems, modifying the system dynamics using control is a more natural perspective. The central philosophy of the paper is to design a control scheme to achieve a particular objective optimally (if possible).
\section*{Appendix: The P\&I Approach}\label{Appendix 2}
Given a single input system
\begin{align}\label{feedbackoriginalsystem}
\begin{split}
\dt{\mathrm{x}}=\mathrm{\mathbb{F}(\mathrm{x}, \lambda)}\hspace{1.5cm}
\dt{\lambda}=\mathrm{u} 
\end{split} 
\end{align}
with $(\mathrm{x}, \lambda)\in (\mathbb{R}^{n-1},\mathbb{R})$ and without any particular structure. 
The target system, the manifold invariance condition, and the implicit manifold condition—i.e., the three essential steps of the classical I\&I approach—are merged in step $(S_1)$. 
The notion of tangent space and passivity theory is invoked to propose the P\&I approach in the following three steps $(S_2-S_4)$ to ensure the \textit{manifold attractivity}.
 
$(S_1)$ \textbf{Construction of the implicit manifold:}
The target dynamics $\dt{\eta}=\beta(\eta)$  with $\mathrm{x}=\eta$ is defined such that the subsystem $\dt{\mathrm{x}}=\mathrm{\mathbb{F}(\mathrm{x}, \varphi (\mathrm{\mathrm{x}}))}$ for $\mathrm{C}^{\infty}$  mapping $\varphi(\mathrm{x}):\mathbb{R}^{\mathrm{n}}\rightarrow \mathbb{R}$ has a GES/GAS equilibrium at the origin by considering the relationship $\mathrm{\lambda}=\varphi (\mathrm{\mathrm{x}})$. This defines the implicit manifold  $\Psi(\mathrm{x}, \lambda )=\mathrm{\lambda}-\varphi (\mathrm{\mathrm{x}})=0$, the implicit manifold $\mathbb{M}=\left \{(\mathrm{x},\lambda) \in \mathbb{R}^{\mathrm{n-1}}\times\mathbb{R} | \Psi(\mathrm{x}, \lambda):=\lambda-\varphi(\mathrm{x})=0 \right \}$, and $\pi(\eta)=\mathrm{col}(\eta, \varphi(\mathrm{\eta}))$.

{$(S_2)$}  \textbf{Tangent space structure for control systems:}
Consider an n-dimensional manifold $\mathbb{M}$ with tangent bundle $\mathbb{T}_{\mathbb{M}}$, such that all $\mathrm{p}\in {\mathbb{M}}$, $\mathbb{T}_{\mathrm{p}}{\mathbb{M}}$ has the following structure
\begin{equation}
    \mathbb{T}_{\mathrm{p}}{\mathbb{M}} = \mathbb{H}_{\mathrm{p}} \oplus  \mathbb{V}_{\mathrm{p}}: \hspace{0.3cm}  \mathbb{H}_{\mathrm{p}} \cap   \mathbb{V}_{\mathrm{p}}=0
\end{equation}
where $\mathbb{H}_{\mathrm{p}}$ is the horizontal space and $\mathbb{V}_{\mathrm{p}}$ is the vertical space.  
\begin{equation}\label{TpMm}
  \mathrm{Then} \hspace{0.9cm} \mathbb{T}_{\mathrm{p}}{\mathbb{M}}= \mathbb{H}_{\mathrm{p}} \oplus  \mathbb{V}_{\mathrm{p}}=(\dt{\mathrm{x}}, 0)\oplus (0, \dt{\lambda})=(\dt{\mathrm{x}},\dt{\lambda})
\end{equation} 
is written for given system (\ref{feedbackoriginalsystem}) at any point $\mathrm{p}\in\mathbb{M}$. With the implicit manifold $\Psi (\mathrm{x} ,\lambda )$ obtained in $S_1$, the normal vector direction is given by $\triangledown \Psi (\mathrm{x} ,\lambda )$. Thus, a PR metric $\mathrm{R}$ on space $\mathbb{T}_{\mathrm{p}}{\mathbb{M}}$ can be defined as
\begin{small}
\begin{align}\label{metric}
    \begin{split}
        \mathrm{R}&=\triangledown \Psi (\mathrm{x} ,\lambda )^{\mathrm{T}}\triangledown \Psi (\mathrm{x} ,\lambda )\\&=\begin{bmatrix}
\left ( \frac{\partial \varphi}{\partial \mathrm{x}}\right )^{\mathrm{T}}\left ( \frac{\partial \varphi}{\partial \mathrm{x}}\right ) & \left ( -\frac{\partial \varphi}{\partial \mathrm{x}}\right )^{\mathrm{T}}\\ 
-\left ( \frac{\partial \varphi}{\partial \mathrm{x}}\right ) & \mathrm{I}
\end{bmatrix}=\begin{bmatrix}
\mathrm{\mathrm{m}}_{11} & \mathrm{\mathrm{m}}_{12}\\ 
 \mathrm{\mathrm{m}}_{21}&\mathrm{\mathrm{m}}_{22}
\end{bmatrix}.
    \end{split}
\end{align}
\end{small}
which is intuitively a natural choice.
\begin{remark}
To obtain the passive output, the metric $\chi$ is replaced with semi-Riemannian metric $\mathrm{R}$ as a natural choice. (Details and proof can be found in \cite{nayyer2022towards})
\end{remark}
For $(\mathbb{M}, \mathrm{R})$, the splitting is visualized as follows:
\begin{small}
\begin{equation}
     (\dt{\mathrm{x}},\dt{\lambda})=\left (\dt{\mathrm{x}}, -\mathrm{m}_{22}^{-1}\mathrm{m}_{21} \dt{\mathrm{x}} \right )\oplus \left (0, \dt{\lambda}+\mathrm{m}_{22}^{-1}\mathrm{m}_{21} \dt{\mathrm{x}} \right ) =\widetilde{\mathbb{H}}_{\mathrm{p}} \oplus \widetilde{\mathbb{V}}_{\mathrm{p}}
\end{equation}
\end{small}

As $\dt{\lambda}$  is along the vertical direction, the passive output is chosen as a component $\dt{\lambda}+\mathrm{m}_{22}^{-1}\mathrm{m}_{21} \dt{\mathrm{x}}$ which is in the same direction or parallel to $\dt{\lambda}$ \cite{nayyer2022towards}.  Roughly speaking, the idea is to bring the component $\dt{\lambda}+\mathrm{m}_{22}^{-1}\mathrm{m}_{21}$ of $\widetilde{\mathbb{H}}_{\mathrm{p}}$ to the component $\mathrm{m}_{22}^{-1}\mathrm{m}_{21} \dt{\mathrm{x}}$ of $\widetilde{\mathbb{V}}_{\mathrm{p}}$.
\begin{figure}[ht!]
    \centering
    \includegraphics[width=\linewidth]{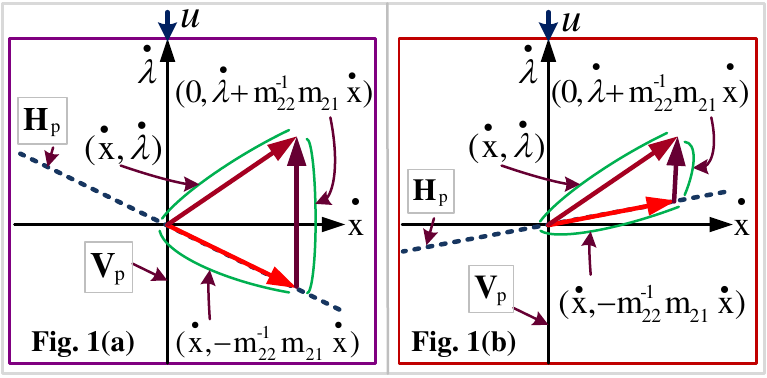}
    \caption{Geometrical interpretation: vertical vector $\mathbb{V}_{\mathrm{p}}$ is along the fiber direction and $\mathbb{H}_{\mathrm{p}} \oplus  \mathbb{V}_{\mathrm{p}}=\mathbb{T}_{\mathrm{p}}{\mathbb{M}}$. }
    \label{Geometrical interpretation}
\end{figure}
The geometrical interpretation for the splitting tangent vector is shown in Fig. \ref{Geometrical interpretation}. 

\textbf{$(S_3)$  Passive output:}
The component of $\mathrm{u}$ tangent vector along $\dt{\lambda}$  is used to define the \textit{passive output} 
\begin{equation}\label{GammaPassive}
    \mathrm{y}=\mathrm{y}_1+\mathrm{y}_2=\int_{0}^{t}\dt{\lambda}\hspace{0.08cm}\mathrm{dt}
+{\int_{0}^{t}(\mathrm{m}_{22}^{-1}\mathrm{m}_{21}\dt{\mathrm{x}})\hspace{0.08cm}\mathrm{dt}}
\end{equation}
with the help of the passivity theory. Here, $\mathrm{y}_1={\int_{0}^{t}\dt{\lambda}\hspace{0.08cm}\mathrm{dt}}$ and $\mathrm{y}_2={\int_{0}^{t}(\mathrm{m}_{22}^{-1}\mathrm{m}_{21}\dt{\mathrm{x}})\hspace{0.08cm}\mathrm{dt}}$ are defined. If $\mathrm{m}_{22}^{-1}\mathrm{m}_{21}$ is a constant then $\mathrm{y}=({\lambda}+\mathrm{m}_{22}^{-1}\mathrm{m}_{21}{\mathrm{x}})$ is defined. If $\mathrm{m}_{22}^{-1}\mathrm{m}_{21}$ is a function of $\mathrm{x}$ then, it can be written as the gradient of any function $\mathrm{q(x)}$ i.e.,  $\mathrm{m}_{22}^{-1}\mathrm{m}_{21}\mathrm{(x)}=\triangledown \mathrm{q(x)}$. Then
\begin{equation}\label{Gamma}
    \mathrm{y}=\int_{0}^{t}(\dt{\lambda}+\triangledown \mathrm{q(x)} \dt{\mathrm{x}})\mathrm{dt}=(\lambda+\mathrm{q(x)})
\end{equation}
\begin{remark}
The condition  $\mathrm{m}_{22}^{-1}\mathrm{m}_{21}\mathrm{(x)}=\triangledown \mathrm{q(x)}$ is related to the condition of integrability and integrable connection in differential geometry. 
\end{remark}

$(S_4)$ \textbf{Storage function:}
With $\mathrm{y}$, the candidate Lyapunov function $\mathbb{S}(\mathrm{x}, \lambda)$ (i.e., storage function) is defined as
\begin{equation}\label{Storage function}
  \mathbb{S}(\mathrm{x}, \lambda)=\frac{1}{2}{\mathrm{y}}^2=\frac{1}{2}(\lambda+\mathrm{q(x)})^2.
\end{equation}
The convergence of the off-the-manifold dynamics to the implicit manifold at an exponential rate $\alpha$ is accompanied by selecting the condition
\begin{equation}\label{exponential laypunov}
    \dt{\mathbb{S}}\leq -\alpha\mathbb{S}.
\end{equation}
One can use the condition (\ref{exponential laypunov}) along with the storage function (\ref{Storage function}) and passive output (\ref{Gamma}) to get
\begin{equation}\label{exponential laypunov solution}
    (\lambda+\mathrm{q(x)})(\dt{\lambda}+\frac{\partial \mathrm{q(x)}}{\partial \mathrm{x}}\dt{\mathrm{x}})\leq -\frac{\alpha}{2}(\lambda+\mathrm{q(x)})^2
\end{equation}
\begin{equation}\label{sjksdkc}
 \Rightarrow (\dt{\lambda}+\frac{\partial \mathrm{q(x)}}{\partial \mathrm{x}}\dt{\mathrm{x}})+\frac{\alpha}{2}(\lambda+\mathrm{q(x)})=0
\end{equation}
The equation (\ref{sjksdkc}) is modified by substituting (\ref{feedbackoriginalsystem}) in terms of the final control law as
\begin{equation}\label{final control law}
  \mathrm{u}=-\frac{\alpha}{2}\lambda-\frac{\alpha}{2}\mathrm{q(x)}-\frac{\partial\mathrm{q(x)}}{\partial \mathrm{x}}\mathrm{f(\mathrm{x}, \lambda)}.
\end{equation}
The above-defined control law ensures the GAS equilibrium point of the system to zero/origin.


\bibliography{Bibliography}
\bibliographystyle{IEEEtran}


\end{document}